# An improved bound on the Minkowski dimension of Besicovitch sets in $\mathbb{R}^3$

By Nets Hawk Katz, Izabella Łaba, and Terence Tao

*In memory of Tom Wolff* (1954–2000)

### Abstract

A Besicovitch set is a set which contains a unit line segment in any direction. It is known that the Minkowski and Hausdorff dimensions of such a set must be greater than or equal to $5/2$ in $\mathbb{R}^3$. In this paper we show that the Minkowski dimension must in fact be greater than $5/2 + \varepsilon$ for some absolute constant $\varepsilon > 0$. One observation arising from the argument is that Besicovitch sets of near-minimal dimension have to satisfy certain strong properties, which we call "stickiness," "planiness," and "graininess."

The purpose of this paper is to improve upon the known bounds for the Minkowski dimension of Besicovitch sets in three dimensions. As a by-product of the argument we obtain some strong conclusions on the structure of Besicovitch sets with almost-minimal Minkowski dimension.

*Definition* 0.1. A *Besicovitch set* (or "Kakeya set") $E \subset \mathbb{R}^n$ is a set which contains a unit line segment in every direction.

Informally, the Kakeya conjecture states that all Besicovitch sets in $\mathbb{R}^n$ have full dimension; this conjecture has been verified for $n = 2$ but is open otherwise. For the purposes of this paper we shall restrict ourselves to the Minkowski dimension, which we now define.

*Definition* 0.2. If $E$ is in $\mathbb{R}^n$, we define the $\delta$-entropy $\mathcal{E}_\delta(E)$ of $E$ to be the cardinality of the largest $\delta$-separated subset of $E$, and $N_\delta(E)$ to be the $\delta$-neighbourhood of $E$.

*Definition* 0.3. For any set $E \subset \mathbb{R}^n$, the (upper) *Minkowski dimension* $\overline{\dim}(E)$ is defined as

$$\overline{\dim}(E) = \limsup_{\delta \to 0} \log_{1/\delta} \mathcal{E}_\delta(E) = n - \liminf_{\delta \to 0} \log_\delta |N_\delta(E)|.$$

In $\mathbb{R}^n$, Wolff [17] showed the estimate

(1) $$\overline{\dim}(E) \geq \frac{1}{2}n + 1,$$

while Bourgain [4] has shown

$$\overline{\dim}(E) \geq \frac{13}{25}n + \frac{12}{25}.$$

The latter result has recently been improved in [8] to

(2) $$\overline{\dim}(E) \geq \frac{4}{7}n + \frac{3}{7}.$$

For further results, generalizations, and applications see [19].

When $n = 3$ Wolff's bound is superior, and thus the best previous result on the three-dimensional problem was $\overline{\dim}(E) \geq 5/2$. By combining the ideas of Wolff and Bourgain with some observations on the structure of (hypothetical) extremal counterexamples to the Kakeya problem, we have obtained the following improvement, which is the main result of this paper.

THEOREM 0.4. *There exists an $\varepsilon > 0$ such that $\overline{\dim}(E) \geq 5/2 + \varepsilon$ for all Besicovitch sets $E$ in $\mathbb{R}^3$.*

While the epsilon in this theorem could in principle be computed, we have not tried to optimize our arguments in order to produce an efficient value for $\varepsilon$. The argument in this paper certainly works for $\varepsilon = 10^{-10}$, but this is definitely far from best possible.

Broadly speaking, the argument is a proof by contradiction. A hypothetical counterexample to the theorem is assumed to exist for some small $\varepsilon$. By refining the collection of tubes slightly, and at one point passing from scale $\delta$ to scale $\rho = \sqrt{\delta}$, one can impose a surprisingly large amount of structure on "most" of the Besicovitch set. Eventually there will be enough structure that one can apply the techniques of Bourgain[4] efficiently and obtain a contradiction. (Of course, if one applies these techniques directly then one would only obtain (2), which is inferior for $n = 3$.)

By the term "most" used in the previous paragraph, we roughly mean that the portion of the Besicovitch set for which our structural assumptions fail only occupies an extremely small fraction of the entire set; we will make this notion precise in Section 5. We shall need this very strong control on the exceptional set, as there is a key stage in the argument in which we need to find an arithmetic progressions of length three in the nonexceptional portion of the Besicovitch set. A discussion of the difficulties of this approach when one only knows that a small portion of each tube is "good" can be found in Bourgain [4]. This means that we will not use methods such as pigeonholing to obtain structural assumptions on our set, as these types of methods usually only give a nonexceptional set which is about $\log(1/\delta)^{-1}$ of the full Besicovitch set.

It may well be that one can use further ideas in [4], such as using triples of points whose *reciprocals* are in arithmetic progression, in order to circumvent this restriction. However, there is an additional obstruction preventing us from obtaining improvements to Theorem 0.4 such as a Hausdorff dimension or maximal function result, as in [17] or [4]. Namely, our argument crucially requires control of the entropy of the Besicovitch set not only at scale $\delta$, but also at many intermediate scales between $\delta$ and 1. In particular, the scale $\rho = \sqrt{\delta}$ plays a key role. Such control is readily available in the case when the upper Minkowski dimension is assumed to be small, but not in the other cases just discussed.

We will derive structural properties on our hypothetical low-dimensional Besicovitch set in the following order. Firstly, we follow an observation in Wolff [18] and observe that the Besicovitch set must be "sticky," which roughly states that the map from directions to line segments in the Besicovitch set is almost Lipschitz. To make this observation rigorous we require the X-ray estimate in [18], and also rely crucially on the fact mentioned earlier, that we have control of the Besicovitch set at multiple scales. We will achieve stickiness in Section 3.

Once we have obtained stickiness, it is a fairly routine matter to show that the Besicovitch set must behave in a self-similar fashion. For instance, if one takes a $\rho$-tube centered around one of the line segments in the Besicovitch set, and dilates it by $\rho^{-1}$ around its axis, one should obtain a new Besicovitch set with almost identical properties. We will obtain quantitative versions of these heuristics in Section 6, after some preliminaries in Sections 4 and 5.

We then combine these self-similarity properties with the following geometric heuristic: if for $i = 1, 2, 3$, we have a vector $v_i$ and a family of tubes $\mathbb{T}_i$ which all approximately point in the direction $v_i$, then the triple intersection

$$(3) \qquad \bigcap_{i=1}^{3}\left(\bigcup_{T\in\mathbb{T}_i} T\right)$$

will be fairly small unless $v_1, v_2, v_3$ are almost co-planar. As a consequence we be able to conclude a remarkable structural property on the Besicovitch set, which we call "planiness." Roughly speaking, it asserts that for most points $x$ in the Besicovitch set, most of the line segments passing through $x$ lie on a plane $\pi(x)$, or on the union of a small number of planes. At first we shall only derive this property at scale $\rho$, for reasons which shall become clear, but by changing scale we may easily impose this property at scale $\delta$ as well.

One can analyze the derivation of planiness further, and obtain an important additional property which we call "graininess." Roughly speaking, this asserts that the intersection of the Besicovitch set with any $\rho$-cube will, when studied at scale $\delta$, look like a union of $\delta \times \rho \times \rho$ boxes which are parallel to the plane $\pi(x)$ mentioned earlier. This property is obtained by repeating the derivation of planiness, but with the additional observation that even if $v_1, v_2, v_3$ lie in a common plane $\pi(x)$, the set (3) will still be small unless the sets $\bigcup_{T \in \mathbb{T}_i} T$ are essentially of the form just described, assuming that the angles between $v_1$, $v_2$, and $v_3$ are fairly large.

We will define the properties of planiness and graininess rigorously in Section 8, after some preliminaries in Section 7. The derivation of these properties shall be the most technical part of the paper, requiring an increasingly involved sequence of definitions, but once these properties are obtained, the argument will become technically much simpler (though still somewhat lengthy).

We remark that the arguments up to this point are not necessarily restricted to sets of dimension close to 5/2, although for other dimensions one must find an analogue of Wolff's X-ray estimate [18] to begin the argument. The argument below, however, is only effective near dimension 5/2.

Next, we follow the philosophy of Bourgain [4] and find three $\rho$-cubes $Q_0$, $Q_1$, $Q_2$ in arithmetic progression which each satisfy certain good properties. To do this it is important that all the properties attained up to this point occur on a very large fraction of the set, which unfortunately causes the arguments in previous sections to be somewhat involved. However, once we have obtained the arithmetic progression then one can be far less stringent, and deal with properties that are satisfied fairly sparsely.

We now apply the ideas in [4], which we now pause to recall. Very roughly, the argument in [4] for the Minkowski dimension runs as follows. Let $A$, $B$, $C$ be the intersections of the Besicovitch set with three planes in arithmetic progression. As the Besicovitch set contains many line segments through $A$, $B$, $C$, it follows that there are many pairs $(a, b)$ of points in $A \times B$ whose midpoint is in $C$. In fact, if the dimension of the Besicovitch set is close to $(n + 1)/2$, then a large fraction of $A \times B$ will have this property. Schematically, we may write this property as

$$(A + B) \subset 2C.$$

In particular, since we expect $A$, $B$, $C$ to be of comparable size, we should have

(4) $$|A + B| \approx |A| \approx |B|,$$

where we need to discretize $A$, $B$ at some scale $\delta$ to make sense of the above expressions. From the combinatorial lemmas in [4] relating sums to differences, we should therefore have (in an appropriate sense)

$$|A - B| \ll |A||B|. \tag{5}$$

But this implies that many of the line segments connecting $A$ with $B$ are parallel, which contradicts the definition of a Besicovitch set. Hence we cannot have a Besicovitch set with dimension close to $(n+1)/2$.

Suppose we applied the above ideas to our situation. If $Q_0$, $Q_1$, and $Q_2$ are three $\rho$-cubes in arithmetic progression, the Besicovitch set property should imply that there are many pairs $(a, b)$ of points in $Q_0 \times Q_2$ whose midpoint is in $Q_1$. Unfortunately, if our set has dimension close to $5/2$, this set of pairs of points is very sparse compared to $Q_0 \times Q_2$, and the combinatorial lemmas in [4] are not effective in this context.

However, this can be salvaged by using the planiness and graininess properties of our set, and especially the fact that the planes through $x$ are almost always parallel to the squares through $x$. These very restrictive properties drastically reduce the possible degrees of freedom of the Besicovitch set, and many pairs in $Q_0 \times Q_2$ can be ruled out *a priori* as being of the form above. The end result is that one can find a well-behaved subset $G$ of $Q_0 \times Q_2$ which is determined by the planes and grains, such that the midpoint of $a$ and $b$ is in $Q_1$ for a large fraction of pairs $(a, b)$ in $G$. By applying the lemmas of [4] we find once again that the Besicovitch set contains many tubes which are parallel, if $\varepsilon$ is sufficiently small. This is our desired contradiction, and we are done.

For completeness we also give in an Appendix a sketch of the argument from [4] that we use.

The properties of planiness and graininess may seem strange, but there is a simple example of an object which resembles a $5/2$-dimensional Kakeya set — albeit in $\mathbb{C}^3$ instead of $\mathbb{R}^3$ — and which does obey these properties. Namely, the Heisenberg group

$$\{(z_1, z_2, z_3) \in \mathbb{C}^3 : \mathrm{Im}(z_3) = \mathrm{Im}(z_1 \overline{z_2})\}.$$

has real dimension 5, contains a four-parameter family of lines (some of which, though, are parallel), and satisfies the planiness and graininess properties perfectly. We will discuss this example further in Section 13.

The authors are indebted to Jean Bourgain for explaining his recent work to one of the authors, the referee for helpful comments, and especially to Tom Wolff for his constant encouragement and mathematical generosity. Many of the "new" ideas in this paper were inspired, directly or indirectly, by the (mostly unpublished) heuristics, computations, and insights which Tom Wolff kindly shared with the authors. The first and third authors are supported by NSF grants DMS-9801410 and DMS-9706764 respectively.

## 1. Notation and preliminaries I

Throughout this paper we shall always be working in three dimensions $\mathbb{R}^3$ unless otherwise specified. We use italic letters $x, y, z$ to denote points in $\mathbb{R}^3$, and Roman letters $(\mathrm{x}, \mathrm{y}, \mathrm{z})$ to denote coordinates of points in $\mathbb{R}^3$.

Unless otherwise specified, all integrals will be over $\mathbb{R}^3$ with Lebesgue measure.

In this paper $\delta$ refers to a number such that $0 < \delta \ll 1$, and $\varepsilon$ refers to a fixed number such that $0 < \varepsilon \ll 1$. $d$ will be a number such that $2 < d < 3$; later on we will set $d = 5/2$. The symbol $\rho$ will always denote the quantity $\rho = \sqrt{\delta}$.

We use $C$, $c$ to denote generic positive constants, varying from line to line (unless subscripted), which are independent of $\varepsilon$, $\delta$, but which may depend on the parameter $d$. $C$ will denote the large constants and $c$ will denote the small constants.

We will use $X \lesssim Y$, $Y \gtrsim X$, or $X = O(Y)$ to denote the inequality $X \leq AY$, where $A$ is a positive quantity which may depend on $\varepsilon$. We use $X \gg Y$ to denote the statement $X \geq AY$ for a large constant $A$. We use $X \sim Y$ to denote the statement that $X \lesssim Y$ and $Y \lesssim X$.

We will use $X \lessapprox Y$, $Y \gtrapprox X$, or "$Y$ majorizes $X$" to denote the inequality

$$X \leq A\delta^{-C\varepsilon}Y,$$

where $A$ is a positive quantity which may depend on $\varepsilon$, and $C$ is a quantity which does not depend on $\varepsilon$. We use $X \approx Y$ to denote the statement that $X \lessapprox Y$ and $Y \lessapprox X$. In particular we have $\varepsilon \approx 1$.

If $E$ is a subset of $\mathbb{R}^n$, we use $|E|$ to denote its Lebesgue measure; if $I$ is a finite set, we use $\#I$ to denote its cardinality.

For technical reasons, we will require a nonstandard definition of a $\delta$-tube. Namely, a $\delta$-tube $T$ is a $\delta$-neighbourhood of a line segment whose endpoints $x_0$ and $x_1$ are on the planes $\{(\mathrm{x},\mathrm{y},\mathrm{z}) : \mathrm{z} = 0\}$ and $\{(\mathrm{x},\mathrm{y},\mathrm{z}) : \mathrm{z} = 1\}$ respectively, and whose orientation is within $\frac{1}{10}$ of the vertical. Note that

(6) $$|T_\sigma| \sim \sigma^2$$

for any $\sigma$-tube $T_\sigma$. We define a *direction* to be any quantity of the form $(\mathrm{x}, \mathrm{y}, 1)$ with $|\mathrm{x}|, |\mathrm{y}| \lesssim 1$, and $\mathrm{dir}(T)$ to be the direction $x_1 - x_0$.

If $T$ is a tube, we define $CT$ to be the dilate of $T$ about its axis by a factor $C$. We say that two tubes $T$ and $T'$ are *equivalent* if $T \subset CT'$ and $T' \subset CT$. If $\mathbb{T}$ is a set of tubes, we say that $\mathbb{T}$ *consists of essentially distinct tubes* if for any $T \in \mathbb{T}$ there are at most $O(1)$ tubes $T'$ which are equivalent to $T$.

We use the term $r$-ball to denote a ball of radius $r$, and use $B(x, r)$ to denote the $r$-ball centered at $x$. By a simple covering argument using $\delta$-balls, we see that the quantities $\mathcal{E}_\delta(E)$ and $N_\delta(E)$ defined in Definition 0.2 are related in $\mathbb{R}^3$ by the basic estimate

$$\delta^3 \mathcal{E}_\delta(E) \sim |N_\delta(E)|. \tag{7}$$

If $\mu$ is a function, we use $\mathrm{supp}(\mu)$ to denote the support of $\mu$.

If $1 < d \leq \infty$ is an exponent, we define the dual exponent by $d' = d/(d-1)$, with $\infty' = 1$. The support, $L^1$ norm, $L^{d'}$ norm, and size of a function are all related of course by such standard inequalities as Hölder and Chebyshev. We will rely on these inequalities extremely often, and so we write them down for future reference.

LEMMA 1.1. *For any nonnegative function $\mu$ on a measure space (continuous or discrete), any $1 < d \leq \infty$, and any $\lambda > 0$, we have*

$$\|\mu\|_1 \lessapprox \|\mu\|_{d'} |\mathrm{supp}(\mu)|^{1/d}, \tag{8}$$

$$\|\mu\|_{d'} \gtrapprox \|\mu\|_1 |\mathrm{supp}(\mu)|^{-1/d}, \tag{9}$$

$$\int_{\mu \gtrapprox \lambda} \mu \lessapprox \lambda^{1-d'} \|\mu\|_{d'}^{d'}, \tag{10}$$

$$\int_{\mu \lessapprox \lambda} \mu \lessapprox \lambda |\mathrm{supp}(\mu)|, \tag{11}$$

$$|\{\mu \gtrapprox \lambda\}| \lessapprox \lambda^{-d'} \|\mu\|_{d'}^{d'}. \tag{12}$$

We remark that all the above quantities will automatically be finite in our applications.

Finally, we observe the following trivial uniformity lemma.

LEMMA 1.2. *Let $\mu$ be a nonnegative function on a measure space (continuous or discrete), such that*

$$|\mathrm{supp}(\mu)| \leq A,$$
$$\|\mu\|_\infty \leq B,$$
$$\|\mu\|_1 \gtrapprox AB$$

*for some $A, B > 0$. Then there exists a nonempty set $E \subset \mathrm{supp}(\mu)$ such that*

$$|E| \approx |\mathrm{supp}(\mu)| \approx A,$$
$$\mu \approx B \text{ on } E,$$
$$\|\mu\|_{L^1(E)} \approx AB.$$

*Proof.* Define $\lambda$ by
$$\lambda = (AB)^{-1}\|\mu\|_1;$$
from hypothesis we have $\lambda \approx 1$. We define $E$ to be the set
$$E = \left\{\mu > \frac{1}{2}\lambda B\right\}.$$
From the estimate
$$\lambda AB = \|\mu\|_{L^1(E)} + \|\mu\|_{L^1(E^c)} \leq B|E| + \frac{1}{2}\lambda BA$$
which follows from the hypotheses, we see that $|E| \gtrsim A$. The verification of the remainder of the properties are now routine. □

It is of course possible to make this lemma more precise (e.g. by the pigeonhole principle), but we shall not do so here.

## 2. Kakeya estimates

In this section we summarize the Kakeya and X-ray estimates which we shall need. In the following $\sigma, \theta$ are quantities such that $\delta \leq \sigma \leq \theta \leq 1$.

*Definition* 2.1. If $\mathbb{T}_\sigma$ is a collection of $\sigma$-tubes, we define the *directional multiplicity* $m = m(\mathbb{T}_\sigma)$ to be the largest number of tubes in $\mathbb{T}_\sigma$ whose directions all lie in a cap of radius $\sigma$. If $m \approx 1$, we say that $\mathbb{T}_\sigma$ is *direction-separated*.

*Definition* 2.2. Let $2 < d < 3$. We say that there is an X-*ray estimate at dimension* $d$ if there exist $0 < \alpha, \beta < 1$ for which the following statement holds: For any $\delta$-separated set $\mathcal{E}$ of directions and any collection $\mathbb{T}$ of essentially distinct tubes pointing in directions in $\mathcal{E}$,

$$\left\|\sum_{T \in \mathbb{T}} \chi_T\right\|_{d'} \lesssim \delta^{1-\frac{3}{d}} m^{1-\beta}(\delta^2 \#\mathcal{E})^\alpha, \tag{13}$$

where $m$ is the directional multiplicity of $\mathbb{T}$. If we only assume that (13) holds with $0 < \alpha < 1$ and $\beta = 0$, then we say that there is a *Kakeya estimate at dimension* $d$.

Clearly, an X-ray estimate is stronger than a Kakeya estimate at the same dimension.

THEOREM 2.3 ([17], [18]). *There is a X-ray estimate at dimension* $5/2$.

Indeed, (13) is proven in [17] for $(\alpha, \beta, d) = (7/10, 0, 5/2)$, while in [18] this is improved to $(\alpha, \beta, d) = (7/10, 1/4, 5/2)$. Although these values of $\alpha$ and $\beta$ are sharp for this value of $d$, their exact values are not particularly important for our purposes. We remark that an estimate with $\alpha = 0$ can automatically be improved to an estimate for some positive $\alpha$ thanks to Nikishin-Pisier factorization theory; a discussion of this phenomenon may be found in [1].

We will usually only rely on the following variants of the above estimate:

LEMMA 2.4. *Let $\delta \leq \sigma \leq \theta \ll 1$, and let $\mathbb{T}_\sigma$ be a collection of $\sigma$-tubes whose set of directions all lie in a cap of radius $\theta$. Let $2 < d < 3$ be fixed.*

- *If we have a Kakeya estimate at some dimension d, and if the collection $\mathbb{T}_\sigma$ is direction-separated, then*

$$\Big\| \sum_{T_\sigma \in \mathbb{T}_\sigma} \chi_{T_\sigma} \Big\|_{d'} \lessapprox \sigma^{\frac{d-3}{d}} \theta^{\frac{d+1}{d}}. \tag{14}$$

- *If we have an X-ray estimate at some dimension d, and if $\mathbb{T}_\sigma$ consists of essentially distinct tubes, then*

$$\Big\| \sum_{T_\sigma \in \mathbb{T}_\sigma} \chi_{T_\sigma} \Big\|_{d'} \lessapprox \sigma^{\frac{d-3}{d}} \theta^{\frac{d+1}{d}} m^{1-\beta} \tag{15}$$

*for some $\beta > 0$, where $m$ is the directional multiplicity of $\mathbb{T}_\sigma$.*

*Proof.* We prove only the second claim, as the first follows by setting $m \approx 1$ and $\beta = 0$. By an affine transformation we may assume that $\omega = (0, 0, 1)$. Apply the nonisotropic dilation $(x, y, z) \mapsto (\theta^{-1}x, \theta^{-1}y, z)$. This transforms $\mathbb{T}_\sigma$ to a collection $\mathbb{T}'$ of $\sigma\theta^{-1}$ tubes pointing in a $\sigma\theta^{-1}$-separated set of directions, without significantly affecting the directional multiplicity. From (13) we have

$$\Big\| \sum_{T \in \mathbb{T}'} \chi_{T'} \Big\|_{d'} \lessapprox \sigma^{\frac{d-3}{d}} \theta^{\frac{3-d}{d}} m^{1-\beta}.$$

The claim then follows by undoing the dilation. □

In the specific case $d = 5/2$, $\alpha = 7/10$ we can also obtain the above lemma directly from (13).

## 3. The sticky reduction

In the rest of the paper, $2 < d < 3$ will be a number such that there is an X-ray estimate at dimension $d$. In particular, by the results in [18] we may choose $d = 5/2$.

It is well known that Besicovitch sets must have Hausdorff and Minkowski dimensions $\geq d$. The purpose of this section is to show that one can push this observation a bit further, and conclude that sets whose Minkowski dimension is close to $d$ have a certain "sticky" structure. Our arguments crucially rely on the fact that sets with small Minkowski dimension are under control at several scales simultaneously (in particular, at the scales $\delta$ and $\rho = \sqrt{\delta}$); there does not appear to be any obvious way to apply these heuristics to a set with Hausdorff dimension close to $d$, for instance.

*Definition* 3.1. Let $\mathbb{T}$ be a collection of $\delta$-tubes. We say that $\mathbb{T}$ is *sticky* at scale $\delta$, or just *sticky* for short, if it is direction-separated and there exists a collection $\mathbb{T}_\rho$ of direction-separated $\rho$-tubes and a partition of $\mathbb{T}$ into disjoint sets $\mathbb{T}[T_\rho]$ for $T_\rho \in \mathbb{T}_\rho$ such that

(16) $$T \subset T_\rho \text{ for all } T_\rho \in \mathbb{T}_\rho \text{ and } T \in \mathbb{T}[T_\rho]$$

and we have the cardinality estimates

(17) $$\#\mathbb{T} \approx \delta^{-2},$$

(18) $$\#\mathbb{T}_\rho \approx \delta^{-1},$$

(19) $$\#\mathbb{T}[T_\rho] \approx \delta^{-1} \text{ for all } T_\rho \in \mathbb{T}_\rho.$$

We call $\mathbb{T}_\rho$ the collection of parent tubes of $\mathbb{T}$.

For technical reasons we shall also need an iterated version of stickiness:

*Definition* 3.2. Let $\mathbb{T}$ be a collection of direction-separated $\delta$-tubes of cardinality $\approx \delta^{-2}$. We say that $\mathbb{T}$ is *doubly sticky* if it is sticky at scale $\delta$, and its collection $\mathbb{T}_\rho$ of parent tubes is sticky at scale $\rho$.

PROPOSITION 3.3. *Suppose that there is an X-ray estimate at dimension $d$, and that there exists a Besicovitch set $E$ with $\overline{\dim}(E) < d + \varepsilon$. Then for any sufficiently small $\delta$, there exists a doubly sticky collection $\mathbb{T}$ of tubes at scale $\delta$ with parent collection $\mathbb{T}_\rho$ and grandparent collection $\mathbb{T}_{\delta^{1/4}}$, such that*

(20) $$\Big|\bigcup_{T \in \mathbb{T}} T\Big| \lessapprox \delta^{3-d},$$

(21) $$\Big|\bigcup_{T_\rho \in \mathbb{T}_\rho} T_\rho\Big| \lessapprox \rho^{3-d},$$

(22) $$\Big|\bigcup_{T_{\delta^{1/4}} \in \mathbb{T}_{\delta^{1/4}}} T_{\delta^{1/4}}\Big| \lessapprox \delta^{(3-d)/4}.$$

*Proof.* Let $E$ be a Besicovitch set with Minkowski dimension at most $d+\varepsilon$, and let $\delta \ll 1$ be a fixed. We may assume without loss of generality that $E$ is contained in a fixed ball $B(0, C)$. Then by Definition 0.3 we have

(23) $$|N_\sigma(E)| \lessapprox \sigma^{3-d} \text{ for all } \delta \leq \sigma \leq 1.$$

By taking a $\delta$-separated set of directions oriented within $\frac{1}{10}$ of the vertical, and looking at the associated line segments in $E$, we can find (possibly after rescaling $E$ slightly) a direction-separated set $\mathbb{T}$ of $\delta$-tubes satisfying (17) such that each tube is contained in $N_{C\delta}(E)$. Unfortunately this collection need not be sticky, let alone doubly sticky. To remedy this we shall prune $\mathbb{T}$ of its nonsticky components.

Let $\mathcal{E}$ be a maximal $\rho$-separated set of directions. Call a direction $\omega \in \mathcal{E}$ *sticky* if the tubes $T \in \mathbb{T}$ such that $\mathrm{dir}(T) \in B(\omega, \rho)$ can be covered by $O(\delta^{-C_1 \varepsilon})$ $\rho$-tubes pointing in the direction $\omega$; here $C_1$ is a constant to be chosen later.

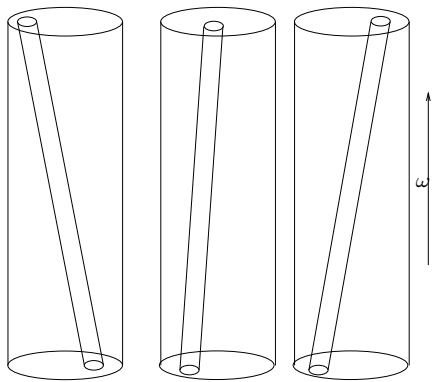

Figure 1. An example of a nonsticky direction $\omega$. The thin tubes are $\delta$-tubes in $N_{C\delta}(E)$, the fat tubes are $\rho$-tubes in $N_{C\rho}(E)$. Note that $N_{C\rho}(E)$ is rather large.

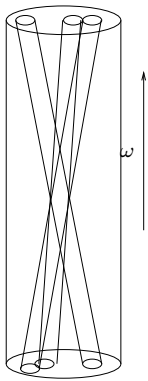

Figure 2. An example of a sticky direction $\omega$.

We now observe that only a small number of directions $\omega$ are nonsticky. More precisely, let $\mathcal{E}_1$ be the subcollection of directions in $\mathcal{E}$ which are nonsticky.

By definition, we can find for each $\omega \in \mathcal{E}_1$ a collection of $\gtrapprox \delta^{-C_1\varepsilon}$ disjoint $\rho$-tubes which are contained in $N_{C\rho}(E)$. Call the union of all these collections $\mathbb{T}'_\rho$. Then by (6)

$$\Big\| \sum_{T \in \mathbb{T}'_\rho} \chi_T \Big\|_1 \sim \rho^2 \# \mathbb{T}'_\rho \gtrapprox \delta^{-C_1\varepsilon}(\rho^2 \# \mathcal{E}_1).$$

From (9) and (23) we thus have

$$\Big\| \sum_{T \in \mathbb{T}'_\rho} \chi_T \Big\|_{d'} \gtrapprox \delta^{-C_1\varepsilon}(\rho^2 \# \mathcal{E}_1)\rho^{-\frac{3-d}{d}}.$$

On the other hand, from (13) we have

$$\Big\| \sum_{T \in \mathbb{T}'_\rho} \chi_T \Big\|_{d'} \lessapprox \rho^{-\frac{3-d}{d}} \delta^{-C_1\varepsilon(1-\beta)} (\rho^2 \# \mathcal{E}_1)^\alpha.$$

Combining the two estimates we obtain

$$\# \mathcal{E}_1 \lessapprox \delta^{-1} \delta^{cC_1\varepsilon},$$

where $c$ depends on $\alpha, \beta$. Since the tubes in $\mathbb{T}$ are direction-separated, we thus see that at most $\lessapprox \delta^{-2}\delta^{cC_1\varepsilon}$ tubes in $\mathbb{T}$ that point within $O(\delta)$ of a nonsticky direction. We may therefore remove these tubes from $\mathbb{T}$ without significantly affecting (17), if $C_1$ was chosen sufficiently large.

The collection $\mathbb{T}$ now has no nonsticky directions. Thus, we may cover the tubes in $\mathbb{T}$ by a family of $\rho$-tubes $\mathbb{T}_\rho$ which are direction-separated. In particular, we have the upper bound in (18). Since $\mathbb{T}$ is direction-separated, each $\mathbb{T}_\rho$ can cover at most $\approx \delta^{-1}$ tubes in $\mathbb{T}$, and so we have the lower bound in (18). Let $\{\mathbb{T}[T_\rho]\}$ be any partition of $\mathbb{T}$ for which (16) always holds. As we have just observed,

$$\# \mathbb{T}[T_\rho] \lessapprox \delta^{-1}$$

for each $T_\rho \in \mathbb{T}_\rho$. On the other hand, from (17) we see that

$$\sum_{T_\rho \in \mathbb{T}_\rho} \# \mathbb{T}[T_\rho] \approx \delta^{-2}.$$

From these estimates, (18), and Lemma 1.2, we can find a subset $\mathbb{T}_\rho'$ of $\mathbb{T}_\rho$ such that $\# \mathbb{T}_\rho' \approx \delta^{-1}$ and (19) holds for all $T_\rho \in \mathbb{T}_\rho'$. If we now replace $\mathbb{T}$ with $\bigcup_{T_\rho \in \mathbb{T}_\rho'} \mathbb{T}[T_\rho]$ and $\mathbb{T}_\rho$ with $\mathbb{T}_\rho'$ we see that we have obtained (19) without significantly affecting any of the other properties just derived.

We now repeat the above procedure but with $\delta$ replaced by $\rho$, and $\mathbb{T}$ replaced by $\mathbb{T}_\rho$. This allows us to refine the collection $\mathbb{T}_\rho$ so that it is also sticky, with an associated collection $T_{\delta^{1/4}}$ of $\delta^{1/4}$-tubes. Of course, to maintain consistency we have to remove the tubes $\mathbb{T}[T_\rho]$ from $\mathbb{T}$ every time we remove a tube $T_\rho$ from $\mathbb{T}_\rho$, but this does not cause any difficulty, and one may verify that all the claims in the proposition hold. Note that the claims (20), (21), and (22) will follow from (23). □

Henceforth we assume that $C_1$ has been set to an absolute constant, and allow all future constants to implicitly depend on $C_1$. We will continue this convention when we choose $C_2$, $C_3$, etc.

The above result established stickiness at the scales $\rho$ and $\delta^{1/4}$. In fact, one could quite easily establish stickiness at every scale $\delta \leq \sigma \leq 1$ from the Minkowski dimension hypothesis, but we shall not need to do so here. We also remark that the above argument requires no special numerology on $d$, and would work perfectly well for $d > 5/2$, providing of course that we had an X-ray estimate at $d$.

For most of the argument we shall not need double stickiness, and derive most of our results just by assuming stickiness. In fact, the only place we shall need double stickiness is in Proposition 9.2, where we need to move the planyness property, which is initially derived at scale $\rho$, to the finer scale of $\delta$ without losing the stickiness property.

## 4. Notation and preliminaries II

In the rest of the paper, $0 < \varepsilon \ll 1$ will be a fixed small number (say $\varepsilon = 10^{-10}$), and $\mathbb{T}$ will be a sticky collection of tubes satisfying (20) and (21). With the exception of Proposition 9.2, we will not use the double-stickiness property, and will not change the value of $\delta$.

For future reference we shall set out some notation and estimates which we shall use frequently. We use $T_\sigma$ to denote a tube of thickness $\sigma$; if $T$ is unsubscripted, we assume it to have thickness $\delta$.

*Definition* 4.1. For any $x \in \mathbb{R}^3$ and $T_\rho \in \mathbb{T}_\rho$, we define the sets $\mathbb{T}(x)$, $\mathbb{T}_\rho(x)$, and $\mathbb{T}[T_\rho](x)$ by

$$\mathbb{T}(x) = \{T \in \mathbb{T} : x \in T\},$$
$$\mathbb{T}_\rho(x) = \{T_\rho \in \mathbb{T}_\rho : x \in T_\rho\},$$
$$\mathbb{T}[T_\rho](x) = \{T \in \mathbb{T}[T_\rho] : x \in T\}.$$

*Definition* 4.2. We define the sets $E_\delta$, $E_\rho$, and $E_\delta[T_\rho]$ for all $T_\rho \in \mathbb{T}_\rho$ by

$$E_\delta = \bigcup_{T \in \mathbb{T}} T,$$

$$E_\rho = \bigcup_{T_\rho \in \mathbb{T}_\rho} T_\rho,$$

$$E_\delta[T_\rho] = \bigcup_{T \in \mathbb{T}[T_\rho]} T.$$

We similarly define the multiplicity functions $\mu_\delta$, $\mu_\rho$, and $\mu_\delta[T_\rho]$ by

$$\mu_\delta(x) = \sum_{T \in \mathbb{T}} \chi_T(x) = \#\mathbb{T}(x),$$

$$\mu_\rho(x) = \sum_{T_\rho \in \mathbb{T}_\rho} \chi_{T_\rho}(x) = \#\mathbb{T}_\rho(x),$$

$$\mu_\delta[T_\rho](x) = \sum_{T \in \mathbb{T}[T_\rho]} \chi_{T_\rho}(x) = \#\mathbb{T}[T_\rho](x).$$

The size of $E_\delta$ and $E_\rho$ is controlled by (20) and (21), while the $L^{d'}$ norms of $\mu_\delta$, $\mu_\rho$, and $\mu_\delta[T_\rho]$ can be controlled by (14). To complement these bounds we have the following precise control on the $L^1$ norms of $\mu_\delta$, $\mu_s d$, and $\mu_\delta[T_\rho]$:

LEMMA 4.3. *Let $T_\rho$ be an element of $\mathbb{T}_\rho$, and $x_0$ be a point in $T_\rho$. Then we have the $L^1$ estimates*

$$\|\mu_\delta\|_1 \approx 1, \tag{24}$$

$$\|\mu_\rho\|_1 \approx 1, \tag{25}$$

$$\|\mu_\delta[T_\rho]\|_1 \approx \delta, \tag{26}$$

$$\|\mu_\delta[T_\rho]\|_{L^1(B(x_0,C\rho))} \approx \delta^{3/2}. \tag{27}$$

*Proof.* The estimates (24), (25), (26) follow from (17), (18), (19), and (6). The proof of (27) is similar and relies on the geometrical observation (cf. Figure 3)

$$\int_{B(x_0,C\rho)} \chi_T \approx \delta^{5/2} \tag{28}$$

for any $T \in \mathbb{T}[T_\rho]$. □

## 5. Properties which occur with probability close to 1

In the sequel we shall frequently need a quantitative version of the statement "The property $P(x)$ implies the property $Q(x)$ with probability close to 1". This motivates

*Definition* 5.1. Let $P(x)$ and $Q(x)$ be logical statements with free parameters $x = (x_1, \ldots, x_n)$, where each of the variables $x_i$ ranges either over a subset of Euclidean space, or over a discrete set. We use

$$\tilde{\forall}_x Q(x) : P(x) \tag{29}$$

to denote the statement that

$$|\{x : Q(x) \text{ holds, but } P(x) \text{ fails}\}| \lessapprox \delta^{c\sqrt{\epsilon}} |\{x : Q(x) \text{ holds}\}| \tag{30}$$

for some absolute constant $c > 0$, where the sets are measured with respect to the measure $dx = \prod_{i=1}^n dx_i$, and $dx_i$ is Lebesgue measure if the $x_i$ range over a subset of Euclidean space, or counting measure if they range over a discrete set.

In practice our variables $x_i$ will either be points in $\mathbb{R}^3$ (and thus endowed with Lebesgue measure), or tubes in $\mathbb{T}$ or $\mathbb{T}_\rho$ (and thus endowed with counting measure). Thus, for instance,

$$\tilde{\forall}_{T,x} T \in \mathbb{T}, x \in T : P(x, T)$$

denotes the statement that

$$\sum_{T \in \mathbb{T}} |\{x \in T : P(x, T) \text{ fails}\}| \lessapprox \delta^{c\sqrt{\epsilon}} \sum_{T \in \mathbb{T}} |T|.$$

The right-hand side of (30) will always be automatically finite in our applications. Note that (29) vacuously holds if $P(x)$ is never satisfied.

Of course, (30) is trivial if $c = 0$. The reason why we choose the factor $\delta^{c\sqrt{\epsilon}}$ is that it is much smaller than any quantity of the form $\delta^{C\varepsilon}$, but much larger than anything of the form $\delta^c$. In particular, (30) implies

$$|\{x : Q(x), P(x)\}| \sim |\{x : Q(x)\}|. \tag{31}$$

Here and in the rest of the paper, the expression "$Q(x), P(x)$" is an abbreviation for "$Q(x)$ and $P(x)$ both hold."

The symbol $\tilde{\forall}_x$ should be read as "for most $x$ such that," where "most" means that the event occurs with probability very close to 1. Observe that the meaning of $\tilde{\forall}_x$ does not depend on the choice of scale ($\delta$, $\rho$, or $\delta^{1/4}$), except that the $c$ in (30) may change value.

We now develop some technical machinery to manipulate expressions of the form (29). This machinery would all be trivial if $\tilde{\forall}$ were replaced by $\forall$, and are not particularly difficult to prove. We first observe some trivial lemmas.

LEMMA 5.2. *Suppose that $m \approx 1$ is an integer and $Q$ and $P_1, \ldots, P_m$ are properties depending on some free parameters $x_1, \ldots, x_n$, such that*
$$\tilde{\forall}_x Q(x) : P_i(x) \text{ holds for all } i = 1, \ldots, m$$
*where the implicit constants are independent of $i$. Then we have*
$$\tilde{\forall}_x Q(x) : P_1(x), \ldots, P_m(x).$$

*Proof.* Apply (30) for each $P_i(x)$ and sum in $i$. □

Let $Q(x) \implies P(x)$ denote the statement "$Q(x)$ fails, or $P(x)$ holds".

COROLLARY 5.3 (Modus Ponens). *If $R(x)$, $Q(x)$, $P(x)$ are properties such that*
$$\tilde{\forall}_x R(x) : Q(x) \text{ and } \tilde{\forall}_x R(x) : (Q(x) \implies P(x))$$
*hold, then we have*
$$\tilde{\forall}_x R(x) : P(x).$$

LEMMA 5.4. *Let $P$, $Q$, $R$ be properties with free parameters $x = (x_1, \ldots, x_n)$ such that*
$$\tilde{\forall}_x R(x) : Q(x) \text{ and } \tilde{\forall}_x R(x), Q(x) : P(x)$$
*hold. Then we have*
$$\tilde{\forall}_x R(x) : P(x), Q(x).$$

*Proof.* Applying (31) to the first hypothesis, and (30) to the second, we obtain
$$\tilde{\forall}_x R(x) : (Q(x) \implies P(x)),$$
and the claim follows from Corollary 5.3. □

LEMMA 5.5. *Let $P(x)$, $Q(x,y)$, $R(x,y)$ be properties with free parameters $x = (x_1, \ldots, x_n)$ and $y = (y_1, \ldots, y_m)$ such that*
$$\tilde{\forall}_y Q(x,y) : R(x,y)$$
*uniformly for all $x$ satisfying $P(x)$. Then we have $\tilde{\forall}_{x,y} P(x), Q(x,y) : R(x,y)$.*

*Proof.* Apply (30) to the hypothesis, and integrate over all $x$ that satisfy $P(x)$. □

Unlike the case with the more familiar $\forall$ quantifier, some care must be taken with $\tilde{\forall}$ when adding or removing dummy variables if this significantly changes the underlying measure. For instance,
$$\tilde{\forall}_x \; x \in \bigcup_{T \in \mathbb{T}} T : P(x)$$

is not necessarily equivalent to

$$\tilde{\forall}_{T,x} T \in \mathbb{T}, x \in T : P(x)$$

because the underlying measures are quite different. In the former case all points $x$ in $\bigcup_{T \in \mathbb{T}}$ have equal weight, whereas in the latter case the weight of a point $x$ is proportional to the multiplicity $\mu_\delta(x)$.

On the other hand, it is legitimate to add or remove dummy variables when the multiplicity is approximately constant. More precisely:

LEMMA 5.6. *Let $Q_1(x)$, $Q_2(x,y)$ be properties such that*

(32) $$|\{y : Q_2(x,y)\}| \approx M$$

*whenever $Q_1(x)$ holds for some quantity $M$ independent of $x$. Then for any property $P(x)$, the statements*

$$\tilde{\forall}_x Q_1(x) : P(x)$$

*and*

$$\tilde{\forall}_{x,y} Q_1(x), Q_2(x,y) : P(x)$$

*are equivalent.*

*Proof.* From (32) we see that

$$|\{(x,y) : Q_1(x), Q_2(x,y)\}| \approx M|\{x : Q_1(x)\}|$$

and

$$|\{(x,y) : Q_1(x), Q_2(x,y) \text{ hold}, P(x) \text{ fails}\}| \approx M|\{x : Q_1(x) \text{ holds}, P(x) \text{ fails}\}|.$$

The claim then follows by expanding the hypothesis and conclusion using (30). □

COROLLARY 5.7. *Suppose the set $Y$ is partitioned into disjoint subsets $Y[z]$ as $z$ ranges over an index set $Z$. Then for any properties $Q(x,y)$, $P(x,y)$ and any free variables $x$, the statements*

$$\tilde{\forall}_{x,y} y \in Y, Q(x,y) : P(x,y)$$

*and*

$$\tilde{\forall}_{x,y,z} z \in Z, y \in Y[z], Q(x,y) : P(x,y).$$

*are equivalent.*

*Proof.* Apply Lemma 5.6 with $x$, $y$ equal to $(x,y)$, $z$ respectively, $M$ equal to 1, and $Q_2((x,y),z)$ equal to the property that $z \in Z$ and $y \in Y[z]$. □

We will usually apply this corollary with $Y = \mathbb{T}$, $Z = \mathbb{T}_\rho$.

Next, we show how a compound $\tilde{\forall}$ quantifier can be split up into two simpler quantifiers.

LEMMA 5.8. *Suppose that $Q(x,y)$, $P(x,y)$ are properties depending on some free parameters $x = (x_1, \ldots, x_n)$, $y = (y_1, \ldots, y_m)$. Then the statements*

(33) $$\tilde{\forall}_{x,y} Q(x,y) : P(x,y)$$

*and*

(34) $$\tilde{\forall}_{x,y'} Q(x,y') : [\tilde{\forall}_y Q(x,y) : P(x,y)]$$

*are equivalent (up to changes of constants).*

Whenever we have nested $\tilde{\forall}$ statements as in (34), we always assume the implicit constants in the inner $\tilde{\forall}$ to be independent of the variables in the outer $\tilde{\forall}$.

*Proof.* Let $R(x)$ denote the property that

$$\tilde{\forall}_y Q(x,y) : P(x,y)$$

holds.

Assume first that (33) held. If $R(x)$ failed, then

$$|\{y : Q(x,y) \text{ holds}, P(x,y) \text{ fails}\}| \gtrapprox \delta^{c'\sqrt{\epsilon}} |\{y : Q(x,y)\}|$$

for some $c'$ which we shall choose later. Integrating this over all $x$ for which $R(x)$ failed, we obtain

$$|\{(x,y) : Q(x,y) \text{ holds}, P(x,y) \text{ fails}\}| \gtrapprox \delta^{c'\sqrt{\epsilon}} |\{(x,y) : Q(x,y) \text{ holds}, R(x) \text{ fails}\}|.$$

From (33) we therefore have (if $c'$ is chosen sufficiently small)

$$|\{(x,y) : Q(x,y) \text{ holds}, R(x) \text{ fails}\}| \lessapprox \delta^{c\sqrt{\epsilon}} |\{(x,y) : Q(x,y)\}|,$$

and (34) follows by replacing $y$ with $y'$.

Now suppose that (34) held. We need to show that

(35) $$\int |\{y : Q(x,y) \text{ holds}, P(x,y) \text{ fails}\}| \, dx \lessapprox \delta^{c\sqrt{\epsilon}} \int |\{y : Q(x,y)\}| \, dx.$$

The contribution to the left-hand side of (35) when $R(x)$ holds is acceptable, by the definition of $R(x)$. The contribution when $R(x)$ fails is majorized by

$$\int_{R(x) \text{ fails}} |\{y : Q(x,y)\}| = \int |\{y' : Q(x,y') \text{ holds}, R(x) \text{ fails}\}| \, dx$$

and this is acceptable by (34). □

COROLLARY 5.9. *Suppose that $Q_1(x)$, $Q_2(x,y)$, and $P(x,y)$ are properties depending on some free parameters $x = (x_1, \ldots, x_n)$, $y = (y_1, \ldots, y_m)$ which obey (32). Then, the statements*

(36) $$\tilde{\forall}_{x,y} Q_1(x), Q_2(x,y) : P(x,y)$$

*and*

$$\tilde{\forall}_x Q_1(x) : [\tilde{\forall}_y Q_2(x,y) : P(x,y)] \tag{37}$$

*are equivalent (up to changes of constants).*

*Proof.* Assume that (36) held. By Lemma 5.8 we have

$$\tilde{\forall}_{x,y'} Q_1(x), Q_2(x,y') : [\tilde{\forall}_y Q_1(x), Q_2(x,y) : P(x,y)].$$

The second $Q_1(x)$ is redundant. Eliminating the $y'$ variable using Lemma 5.6, one obtains (37). The converse implication follows by reversing the above steps. □

We now apply the above machinery to our specific setting, in which we have a sticky collection of tubes.

LEMMA 5.10. *Let $\mathbb{T}$ be a sticky collection of tubes, and let $P(y, T_\rho, T)$ be a property. Then the statements*

$$\tilde{\forall}_{T_\rho, T, y} T_\rho \in \mathbb{T}_\rho, T \in \mathbb{T}[T_\rho], y \in T : P(y, T_\rho, T), \tag{38}$$

$$\tilde{\forall}_{T_\rho, T, y, x} T_\rho \in \mathbb{T}_\rho, T \in \mathbb{T}[T_\rho], y \in T, x \in T_\rho \cap B(y, C\rho) : P(y, T_\rho, T), \tag{39}$$

*and*

$$\tilde{\forall}_{T_\rho, x} T_\rho \in \mathbb{T}_\rho, x \in T_\rho : [\tilde{\forall}_{T,y} T \in \mathbb{T}[T_\rho], y \in T \cap B(x, C\rho) : P(y, T_\rho, T)] \tag{40}$$

*are equivalent (up to changes of constants).*

*Proof.* The equivalence of (38) and (39) follows from Lemma 5.6, since (32) follows from the trivial estimate

$$|T_\rho \cap B(y, C\rho)| \approx \rho^3 \text{ whenever } y \in T \in \mathbb{T}[T_\rho].$$

But if we rearrange (39) as

$$\tilde{\forall}_{T_\rho, x, T, y} T_\rho \in \mathbb{T}_\rho, x \in T_\rho, T \in \mathbb{T}[T_\rho], y \in T \cap B(x, C\rho) : P(y, T_\rho, T),$$

then the equivalence of (39) and (40) follows from Corollary 5.9, since (32) follows from (27). □

## 6. Uniformity and self-similarity

In this section we shall investigate how the stickiness hypothesis, combined with the Kakeya estimate at scale $d$, implies certain self-similarity properties of the set $E_\delta$. We will then show in later sections how these properties imply planiness and graininess properties of the set, with these efforts culminating in Corollary 9.5.

The notation will be as in the previous section. Informally, the main result of this section shall be

HEURISTIC 6.1. *The set $E_\delta$ can be mostly covered by about $\rho^{-d}$ $\rho$-balls $B$, such that the intersection of $E_\delta$ with most of these balls has volume about $\delta^{3-d}\rho^d$. Most points $x \in \bigcup_{T \in \mathbb{T}} T$ are contained in about $\delta^{-(3-d)}$ tubes $T \in \mathbb{T}$, which in turn are mostly contained in about $\rho^{-(3-d)}$ families $\mathbb{T}[T_\rho]$ in such a way that each family contains $\rho^{-(3-d)}$ of the tubes. Also, the set $E_\delta$ is not "lumpy" at any scale greater than $\delta$, in the sense that*

$$|E_\delta \cap B_\sigma| \ll |B_\sigma| \tag{41}$$

*for all $\sigma \gg \delta$ and all $\sigma$-balls $B_\sigma$. Finally, the angle subtended by most pairs of intersecting tubes is $\approx 1$.*

Note that the numerology in this heuristic is consistent with (20), (21), (14), and Lemma 4.3. Moreover, it is essentially the only such numerology which is consistent with these estimates.

The main purpose of this section will be to prove a rigorous version of Heuristic 6.1; this will be done in Proposition 6.6. Of course, in order to interpret the term "most" in the above heuristic, the $\tilde{\forall}$ machinery in the previous section will be used heavily.

We begin with a more accurate version of (41), which may be thought of as a dual to the Minkowski dimension estimate (23).

PROPOSITION 6.2. *For each $\delta \leq \sigma \ll 1$ and $x \in \mathbb{R}^3$, consider the statement*

$$|E_\delta \cap B(x, C\sigma)| \lessapprox \delta^{-\sqrt{\epsilon}} \delta^{3-d} \sigma^d. \tag{42}$$

*Then, if the constants in (42) are chosen appropriately, we have*

$$\tilde{\forall}_{T,x} \, T \in \mathbb{T}, x \in T : (42) \text{ holds for all } \delta \leq \sigma \ll 1.$$

*Proof.* If (42) holds for $\sigma$, and $\sigma' \approx \sigma$, then (42) holds for $\sigma'$ (with slightly worse constants). Thus in order to verify (42) for all $\delta \leq \sigma \ll 1$, it suffices to verify (42) for all $\sigma$ of the form $\sigma = \delta^{k\varepsilon}$, where $k = 0, \ldots, \varepsilon^{-1}$. Thus it suffices

to show that

$$\tilde{\forall}_{T,x} T \in \mathbb{T}, x \in T : \text{(42) holds for all } \sigma = \delta^{k\varepsilon}, k = 0, \ldots, \varepsilon^{-1}.$$

Since $\varepsilon^{-1} \approx 1$, it thus suffices by Lemma 5.2 to show that

$$\tilde{\forall}_{T,x} T \in \mathbb{T}, x \in T : \text{(42) holds for } \sigma$$

uniformly in $\sigma$.

Fix $\sigma$, and let $X$ denote the set

$$X = \{x \in \mathbb{R}^3 : \text{(42) fails for } \sigma\}.$$

Then by (30) and (24) it suffices to show

$$\int_X \mu_\delta \lessapprox \delta^{c\sqrt{\epsilon}}.$$

From (20) and the definition of $X$ we may cover $X$ by a collection $\mathbb{B}$ of $\sigma$-balls with cardinality

(43) $$\#\mathbb{B} \lessapprox \delta^{\sqrt{\epsilon}} \sigma^{-d}.$$

From elementary geometry we have

$$\int_B \chi_T \lessapprox \sigma^{-2} \delta^2 \int_B \chi_{N_{C\sigma}(T)}$$

for all $B \in \mathbb{B}$ and $T \in \mathbb{T}$ (cf. (28)). Summing this in $B$ and $T$ we see that

$$\int_X \mu_\delta \lessapprox \sigma^{-2} \delta^2 \int_{\bigcup_{B \in \mathbb{B}} B} \sum_{T \in \mathbb{T}} \chi_{N_{C\sigma}(T)}.$$

By (8) and (43) this is majorized by

$$\sigma^{-2} \delta^2 [\sigma^3 \delta^{\sqrt{\epsilon}} \sigma^{-d}]^{1/d} \| \sum_{T \in \mathbb{T}} \chi_{N_{C\sigma(T)}} \|_{d'}.$$

The collection of $C\sigma$-tubes $N_{C\sigma(T)}$ can be partitioned into about $\sigma^2 \delta^{-2}$ subcollections, each of which are direction-separated. Thus from the triangle inequality and (14), we can majorize the above by

$$\sigma^{-2} \delta^2 (\sigma^3 \delta^{\sqrt{\epsilon}} \sigma^{-d})^{1/d} \sigma^2 \delta^{-2} \sigma^{\frac{d-3}{d}} = \delta^{\sqrt{\epsilon}/d}$$

as desired. □

We combine this property with two others, which are also related to Heuristic 6.1.

*Definition* 6.3. If $x_0 \in \mathbb{R}^3$ and $T_\rho \in \mathbb{T}_\rho[x_0]$, then we say that $P_1(x_0, T_\rho)$ holds if the three statements

(44) $$|E_\delta[T_\rho] \cap B(x_0, C\rho)| \gtrapprox \delta^{\sqrt{\epsilon}} \delta^{3-d} \rho^d,$$

(45) $\tilde{\forall}_{T,x} T \in \mathbb{T}[T_\rho], x \in T \cap B(x_0, C\rho) : (42)$ holds for all $\delta \leq \sigma \ll 1$,

(46) $\tilde{\forall}_{T,x} T \in \mathbb{T}[T_\rho], x \in T \cap B(x_0, C\rho) : (47)$ holds

hold, where (47) is the estimate

(47) $$\mu_\delta[T_\rho](x) \lessapprox \delta^{-\sqrt{\epsilon}} \rho^{-(3-d)}.$$

PROPOSITION 6.4. *If the constants in Definition* 6.3 *are chosen appropriately, then*

$$\tilde{\forall}_{T_\rho, x_0} T_\rho \in \mathbb{T}_\rho, x_0 \in T_\rho : P_1(x_0, T_\rho).$$

*Proof.* If we rearrange Proposition 6.2 as

$$\tilde{\forall}_{T_\rho, T, x} T_\rho \in \mathbb{T}_\rho, T \in \mathbb{T}[T_\rho], x \in T : (42) \text{ holds for all } \delta \leq \sigma \leq 1$$

and apply Lemma 5.10, we obtain

$$\tilde{\forall}_{T_\rho, x_0} T_\rho \in \mathbb{T}_\rho, x_0 \in T_\rho : (45) \text{ holds}.$$

Next, we show

$$\tilde{\forall}_{T_\rho, x_0} T_\rho \in \mathbb{T}_\rho, x_0 \in T_\rho : (46) \text{ holds}.$$

By Lemma 5.10 again, it suffices to show that

$$\tilde{\forall}_{T_\rho, T, x} T_\rho \in \mathbb{T}_\rho, T \in \mathbb{T}[T_\rho], x \in T : (47) \text{ holds}.$$

By Lemma 5.5 it suffices to show

$$\tilde{\forall}_{T,x} T \in \mathbb{T}[T_\rho], x \in T : (47) \text{ holds}$$

uniformly in $T_\rho$. Fixing $T_\rho$, we rewrite this as

$$\int_{\mu_\delta[T_\rho] \gtrapprox \delta^{-\sqrt{\epsilon}} \rho^{-(3-d)}} \mu_\delta[T_\rho] \lessapprox \delta^{c\sqrt{\epsilon}} \int \mu_\delta[T_\rho].$$

However, from (14) we have

(48) $$\|\mu_\delta[T_\rho]\|_{d'} \lessapprox \delta^{\frac{d-3}{d}} \rho^{\frac{d+1}{d}}.$$

The claim then follows (with $c = d' - 1$) from (10), (26), and some algebra.

We now show that

$$\tilde{\forall}_{x_0} x_0 \in T_\rho : (44) \text{ holds}$$

for all $T_\rho \in \mathbb{T}_\rho$; the proposition then follows from Lemmas 5.5 and 5.2.

Fix $T_\rho$, and let $X$ denote the set

$$X = \{x_0 \in T_\rho : (44) \text{ fails}\}.$$

By (6), we must show that $|X| \lessapprox \delta^{c\sqrt{\epsilon}} \delta$.

If $x_0 \in X$, then from (9), (27), and the failure of (44), we have
$$\|\mu_\delta[T_\rho]\|_{L^{d'}(B(x_0, C\rho))} \gtrsim \delta^{3/2}(\delta^{\sqrt{\epsilon}}\delta^{3-d}\rho^d)^{-1/d}.$$

Comparing this with (48) and applying a covering argument one obtains the entropy estimate
$$\mathcal{E}_\rho(X) \lesssim \left[\delta^{\frac{d-3}{d}}\rho^{\frac{d+1}{d}}\delta^{-3/2}(\delta^{\sqrt{\epsilon}}\delta^{3-d}\rho^d)^{1/d}\right]^{d'}$$

which simplifies to
$$\mathcal{E}_\rho(X) \lesssim \delta^{\sqrt{\epsilon}d'/d}\delta^{-1/2}.$$

The desired bound on $X$ then follows from (7). $\square$

*Definition* 6.5. Let $x_0$ be a point in $\mathbb{R}^3$. We say that $P_2(x_0)$ holds if one has

(49) $$\delta^{C\sqrt{\epsilon}}\rho^{-(3-d)} \lessapprox \#\mathbb{T}_\rho(x_0) \lessapprox \delta^{-C\sqrt{\epsilon}}\rho^{-(3-d)},$$

(50) $$\tilde{\forall}_{T_\rho} T_\rho \in \mathbb{T}_\rho(x_0) \quad : \quad P_1(x_0, T_\rho),$$

(51) $$|E_\delta \cap B(x_0, C\rho)| \lessapprox \delta^{-C\sqrt{\epsilon}}\delta^{3-d}\rho^d, \text{ and}$$

(52) $$\#\{T_\rho \in \mathbb{T}_\rho(x_0) : \operatorname{dir}(T_\rho) \in B(\omega, \theta)\} \lessapprox \delta^{-C\sqrt{\epsilon}}\theta^c\rho^{-(3-d)}$$

for all directions $\omega$ and all $\delta \leq \theta \ll 1$.

Property (52) basically states that the tubes in $\mathbb{T}_\rho(x_0)$ are not clustered in a narrow angular band.

PROPOSITION 6.6. *If the constants in the above definition are chosen appropriately, then we have*

(53) $$\tilde{\forall}_{T_\rho} T_\rho \in \mathbb{T}_\rho, x \in T_\rho : P_2(x).$$

The reader should compare this rigorous proposition with Heuristic 6.1.

*Proof.* By Lemma 5.2 it suffices to handle each of the properties in $P_2(x)$ separately.

We first observe from Proposition 6.4 and Lemma 5.8 that
$$\tilde{\forall}_{T_\rho', x} T_\rho' \in \mathbb{T}_\rho, x \in T_\rho' : \left[\tilde{\forall}_{T_\rho} T_\rho \in \mathbb{T}_\rho, x \in T_\rho : P_1(x, T_\rho)\right].$$

But this is just a rephrasing of
$$\tilde{\forall}_{T_\rho, x} T_\rho \in \mathbb{T}_\rho, x \in T_\rho : (50) \text{ holds}.$$

Next, we show the upper bound in (49):
$$\tilde{\forall}_{T_\rho, x} T_\rho \in \mathbb{T}_\rho, x \in T_\rho : \mu_\rho(x) \lessapprox \delta^{-C\sqrt{\epsilon}}\rho^{-(3-d)}.$$

We rearrange this using (25) as

$$\int_{\mu_\rho \gtrapprox \delta^{-C\sqrt{\epsilon}} \rho^{-(3-d)}} \mu_\rho \lessapprox \delta^{c\sqrt{\epsilon}}.$$

On the other hand, from (14) we have

(54) $$\|\mu_\rho\|_{d'} \lessapprox \rho^{\frac{d-3}{d}},$$

and the claim then follows from (10).

We now turn to the lower bound in (49):

$$\overset{\tilde{\forall}}{T_\rho, x} T_\rho \in \mathbb{T}_\rho, x \in T_\rho : \mu_\rho(x) \gtrapprox \delta^{C\sqrt{\epsilon}} \rho^{-(3-d)}.$$

We may rearrange this using (25) as

(55) $$\int_{\mu_\rho \lessapprox \delta^{C\sqrt{\epsilon}} \rho^{-(3-d)}} \mu_\rho \lessapprox \delta^{c\sqrt{\epsilon}}.$$

The claim then follows from (11), (21), and the fact that $\mu_\rho$ is supported on $E_\rho$.

Next, we show

$$\overset{\tilde{\forall}}{T_\rho, x} T_\rho \in \mathbb{T}_\rho, x \in T_\rho : (52) \text{ holds for all directions } \omega \text{ and all } \delta \leq \theta \ll 1.$$

By repeating the argument in Proposition 6.2, it suffices to obtain

$$\overset{\tilde{\forall}}{T_\rho, x} T_\rho \in \mathbb{T}_\rho, x \in T_\rho : (52) \text{ holds for all directions } \omega$$

uniformly in $\theta$.

Fix $\theta$; by (25), it suffices to show that

$$\int_X \mu_\rho \lessapprox \delta^{c\sqrt{\epsilon}}$$

where $X$ is the set of points in $\mathbb{R}^3$ for which (52) fails for some $\omega$.

By (8) and (54) it suffices to show that

$$\rho^{\frac{d-3}{d}} |X|^{1/d} \lessapprox \delta^{c\sqrt{\epsilon}},$$

or in other words

$$|X| \lessapprox \delta^{c\sqrt{\epsilon}} \rho^{3-d}.$$

If (52) fails for some $\omega$, it also fails (with slightly different constants) for some $\omega$ in a fixed $c\theta$-separated set. Since this set has cardinality $O(\theta^{-2})$, it suffices to show that

$$|\{x \in \mathbb{R}^3 : (52) \text{ fails for } \omega\}| \lessapprox \delta^{c\sqrt{\epsilon}} \rho^{3-d} \theta^2$$

uniformly in $\omega$.

Fix $\omega$. We must show

$$(56) \qquad \left|\left\{\sum_{T_\rho \in \mathbb{T}_\rho'} \chi_{T_\rho} \gtrapprox \delta^{-C\sqrt{\epsilon}} \theta^c \rho^{-(3-d)}\right\}\right| \lessapprox \delta^{c\sqrt{\epsilon}} \rho^{3-d} \theta^2$$

where $\mathbb{T}_\rho'$ is the set of tubes in $T_\rho$ whose directions lie in $B(\omega, \theta)$.

By (12) and (14), the left-hand side of (56) is majorized by

$$[\delta^{-C\sqrt{\epsilon}} \theta^c \rho^{-(3-d)}]^{-d'} [\rho^{1-\frac{3}{d}} \theta^{1+\frac{1}{d}}]^{d'},$$

which we rearrange as

$$\delta^{Cd'\sqrt{\epsilon}} \theta^{\frac{3-d}{d-1} - cd'} \rho^{3-d} \theta^2.$$

Since $d < 3$, (56) follows if $c$ is chosen sufficiently small.

Finally, we address (51). This is easily dealt with because (51) follows from the other properties in $P_2(x_0)$. To see this, first observe from (49) and (50) that $P_1(x_0, T_\rho)$ holds for at least one tube $T_\rho \in \mathbb{T}_\rho$. From (44) and (45) we thus see that there exists an $x \in B(x_0, C\rho)$ such that (42) holds for all $\delta \leq \sigma \ll 1$. Applying this with $\sigma = C\rho$ gives the result. (Alternatively, one can modify the proof of Proposition 6.2 to obtain (51).) $\square$

In the next few sections we will explore some surprising consequences of the property $P_2(x_0)$, namely planiness and graininess.

## 7. Triple intersections

In this section we fix $x_0$ to be a point in $\mathbb{R}^3$ such that $P_2(x_0)$ holds.

Let $T_\rho$ be a tube in $\mathbb{T}_\rho(x_0)$ such that $P_1(x_0, T_\rho)$ holds. By Proposition 6.6 this situation occurs fairly often.

Let $A(x_0, T_\rho)$ denote the set

$$(57) \qquad A(x_0, T_\rho) = E_\delta[T_\rho] \cap B(x_0, C\rho).$$

From (44) we see that

$$|A(x_0, T_\rho)| \gtrapprox \delta^{C\sqrt{\epsilon}} \delta^{3-d} \rho^d.$$

On the other hand, from (51) we see that

$$\left|\bigcup_{T_\rho \in \mathbb{T}_\rho(x_0)} A(x_0, T_\rho)\right| \lessapprox \delta^{-C\sqrt{\epsilon}} \delta^{3-d} \rho^d.$$

Thus we expect a lot of overlap between the $A(x_0, T_\rho)$. In particular, we expect the size of the triple intersection

$$(58) \qquad\qquad A(x_0, T_\rho^1) \cap A(x_0, T_\rho^2) \cap A(x_0, T_\rho^3)$$

to be quite large for many triples $T_\rho^i$, $i = 1, 2, 3$.

However, it turns out that we can get a nontrivial bound on the size of (58) if tubes $T_\rho^i$ are not coplanar. This is basically because (42) forces the projections of (58) along the directions of $T_\rho^i$ to be small. (For a discussion of the relationship between the volume of a set and the areas of its projections, see [11].)

When the tubes $T_\rho^i$ all lie close to a single plane $\pi$, then one cannot directly improve the trivial bound on (58) given by (51). However, we can show that this trivial bound is only attained when the sets $A(x_0, T_\rho^i)$ are "grainy," in the sense that they consist of a union of squares parallel to $\pi$.

More precisely, we have

DEFINITION 7.1. Let $\pi$ be a plane. We define a *square parallel to $\pi$* to be any rectangular box $Q$ of dimensions $\delta \times \rho \times \rho$ whose two long sides are parallel to $\pi$.

LEMMA 7.2. *Let $x_0$ be a point in $\mathbb{R}^3$ such that $P_2(x_0)$ holds, and let $T_\rho^1, T_\rho^2, T_\rho^3$ be three tubes in $\mathbb{T}_\rho(x_0)$. Let $v_i$ be any vectors such that*

$$(59) \qquad\qquad v_i = \rho \operatorname{dir}(T_\rho^i) + O(\delta).$$

*Let $\pi(v_1, v_2)$ denote the plane spanned by $v_1$ and $v_2$, $1 < M \ll \delta^{-1/4}$ be an arbitrary number, and $F$ be an arbitrary subset of $\mathbb{R}^3$.*

- *(Nonplanar tubes). If $\angle v_1, v_2 \geq M^{-1}$, $\angle \pi(v_1, v_2), v_3 \geq \rho M^2$, and (42) holds for all $x \in F$ and $\delta \leq \sigma \leq 1$, then*

$$(60) \qquad |A(x_0, T_\rho^1) \cap A(x_0, T_\rho^2) \cap A(x_0, T_\rho^3) \cap F| \lessapprox \delta^{-C\sqrt{\epsilon}} \delta^{3-d} \rho^d M^{-c}.$$

- *(Planar tubes). If $\angle v_1, v_2 \geq M^{-1}$ and $\angle \pi(v_1, v_2), v_3 \leq \rho M^2$ then*

$$(61)$$
$$|A(x_0, T_\rho^1) \cap A(x_0, T_\rho^2) \cap A(x_0, T_\rho^3) \cap F| \lessapprox \delta^{-C\sqrt{\epsilon}} \delta^{3-d} \rho^d M^C \left( \sup_Q \frac{|F \cap Q|}{|Q|} \right)^c,$$

*where $Q$ ranges over all squares parallel to $\pi(v_1, v_2)$.*

*Here $c > 0$ and $C > 0$ are absolute constants.*

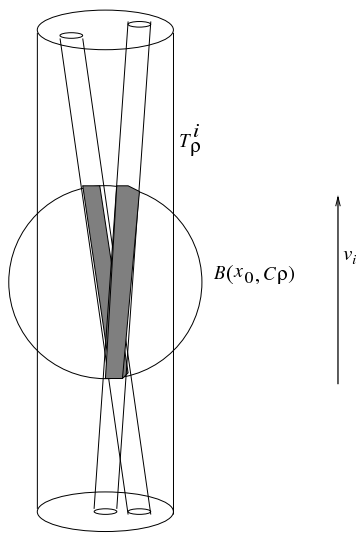

Figure 3. The shaded region is a portion of a typical $A_i$. Note that $A_i$ essentially consists of $\delta \times \rho$ tubes oriented in the direction $v_i$.

*Proof.* Write $A_i = A(x_0, T_\rho^i)$ for short. The key observation is that for every $T \in \mathbb{T}[T_\rho^i]$, the set $T \cap B(x_0, C\rho)$ is essentially constant in the direction $v_i$. More precisely, we have the elementary pointwise estimate (see Figure 3)

(62) $$\chi_{B(x_0,C\rho) \cap T} \lesssim \mathbb{E}_{v_i}(\chi_{B(x_0, 2C\rho) \cap CT})$$

where $E_v$ is the averaging operator

$$\mathbb{E}_v(f(x)) = \int_{|t| \lesssim 1} f(x + vt)\, dt.$$

Summing this in $T$ we obtain

(63) $$\chi_{A_i} \lesssim \mathbb{E}_{v_i}(\chi_{\tilde{A}_i})$$

where

$$\tilde{A}_i = B(x_0, 2C\rho) \cap \bigcup_{T \in \mathbb{T}[T_\rho^i]} CT.$$

It is easy to see that $|\tilde{A}_i| \sim |A_i|$, since $A_i$ is essentially a union of parallel $\delta \times \rho$ tubes and $\tilde{A}_i$ is essentially the union of the dilates of these tubes. In particular, from (51) we have

(64) $$|\tilde{A}_i| \lessapprox \delta^{-C\sqrt{\epsilon}} \delta^{3-d} \rho^d.$$

To utilize (63) we invoke the following estimate, which is a variant of the argument in [11].

LEMMA 7.3. *Let $v_1$, $v_2$, $v_3$ be any three linearly independent vectors, and let $F$ be a subset of $\mathbb{R}^3$. Then for any functions $f_1$, $f_2$, $f_3$,*

$$\int_F \mathbb{E}_{v_1}(f_1)\mathbb{E}_{v_2}(f_2)\mathbb{E}_{v_3}(f_3) \lesssim \sup_P \left(\frac{|F \cap P|}{|P|}\right)^{1/3} \prod_{i=1}^{3} \|f_i\|_3,$$

*where $P$ ranges over all parallelopipeds with edge vectors $v_1$, $v_2$, $v_3$.*

We remark that the factor $\sup_P(\frac{|F \cap P|}{|P|})^{1/3}$ represents a gain over the trivial estimate obtained via Hölder's inequality.

*Proof.* The statement of the lemma is invariant under affine transformations, so we may rescale $v_1$, $v_2$, $v_3$ to be the cardinal directions $e_1$, $e_2$, $e_3$ respectively. It suffices to show that

$$\int_{F \cap P} \mathbb{E}_{v_1}(f_1)\mathbb{E}_{v_2}(f_2)\mathbb{E}_{v_3}(f_3) \lesssim (|F \cap P|)^{1/3} \prod_{i=1}^{3} \|f_i\|_{L^3(CP)},$$

for all unit cubes $P$, since the claim follows by summing over a partition of $\mathbb{R}^3$ and using Hölder's inequality. We may assume that $P$ is centered at the origin, that $F \subset P$ and that $f_i$ are supported on $CP$.

We have the pointwise estimate

$$\mathbb{E}_{v_1}(f_1)(x_1, x_2, x_3) \lesssim \overline{f}_1(x_2, x_3),$$

where $x = (x_1, x_2, x_3)$ and

$$\overline{f}_1(x_2, x_3) = \int_{|x_1| \lesssim 1} f_1(x_1, x_2, x_3) \, dx_1,$$

similarly for cyclic permutations of $1, 2, 3$. Since $\|\overline{f}_i\|_3 \lesssim \|f_i\|_3$ by Young's inequality, we are thus reduced to showing that

$$\int_P \chi_F(x)\overline{f}_1(x_2, x_3)\overline{f}_2(x_3, x_1)\overline{f}_3(x_1, x_2) \, dx_1 dx_2 dx_3 \lesssim |F|^{1/3}\|\overline{f}_1\|_3\|\overline{f}_2\|_3\|\overline{f}_3\|_3.$$

We rewrite this as

$$\int_{|x_1| \lesssim 1} \left( \int_{|x_2|, |x_3| \lesssim 1} \chi_F^{x_1}(x_2, x_3)\overline{f}_1(x_2, x_3)\overline{f}_2^{x_1}(x_3)\overline{f}_3^{x_1}(x_2) \, dx_2 dx_3 \right) dx_1,$$

where $\chi_F^{x_1}(x_2, x_3) = \chi_F(x_1, x_2, x_3)$, etc. By Hölder's inequality this is majorized by

$$\int_{|x_1| \lesssim 1} \|\chi_F^{x_1}\|_3 \|\overline{f}_1\|_3 \|\overline{f}_3^{x_1} \otimes \overline{f}_2^{x_1}\|_3 \, dx_1,$$

where the norms are taken in the $x_2$, $x_3$ variables only. We can simplify this to

$$\|\overline{f}_1\|_3 \int_{|x_1| \lesssim 1} \|\chi_F^{x_1}\|_3 \|\overline{f}_3^{x_1}\|_3 \|\overline{f}_2^{x_1}\|_3 \, dx_1,$$

and the claim follows from applying Hölder in the $x_1$ variable. □

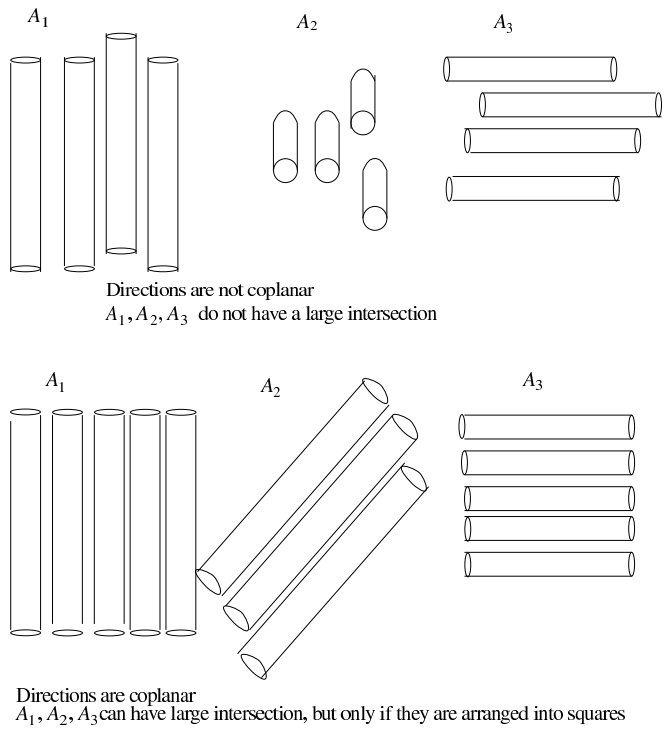

Directions are not coplanar
$A_1, A_2, A_3$ do not have a large intersection

Directions are coplanar
$A_1, A_2, A_3$ can have large intersection, but only if they are arranged into squares

Figure 4. The two cases of Lemma 7.2.

Combining this estimate with (63) and (64) we obtain

$$(65) \qquad |A_1 \cap A_2 \cap A_3 \cap F| \lessapprox \delta^{-C\sqrt{\epsilon}} \delta^{3-d} \rho^d \left( \sup_P \frac{|F \cap P|}{|P|} \right)^{1/3}$$

where $P$ ranges over all parallelopipeds with edge vectors $v_1$, $v_2$, $v_3$.

We now prove the two claims of Lemma 7.2.

*Proof of* (60). Since $A_i \subset E_\delta$, we may assume without loss of generality that $F \subset E_\delta$. Since we are assuming (42) to hold on $F$, we thus see that

$$\frac{|F \cap B|}{|B|} \lessapprox \delta^{-C\sqrt{\epsilon}} \left( \frac{\sigma}{\delta} \right)^{-(3-d)}$$

for any ball of radius $\sigma$.

Let $P$ be any parallelopiped with edge vectors $v_1$, $v_2$, $v_3$. From the hypotheses on $v_1$, $v_2$, $v_3$ and some elementary geometry we see that

$$|P| \sim |N_{M\delta}(P)|,$$

so that $P$ can be efficiently covered by a finitely overlapping collection of $M\delta$-balls. From a simple covering argument we thus obtain

$$\frac{|F \cap P|}{|P|} \lessapprox \delta^{-C\sqrt{\epsilon}} M^{-(3-d)}$$

uniformly in $P$, and (60) follows from (65).

*Proof of* (61). By perturbing $v_i$ if necessary we may assume that $\angle \pi(v_1, v_2), v_3 \gtrsim \delta$. The claim then follows from (65) and the geometric observation that, if $P$ is any parallelopiped with edge vectors $v_1$, $v_2$, $v_3$, then $P$ can be covered by $O(M^2)$ squares $Q$ parallel to $\pi(v_1, v_2)$, and $|P| \gtrsim M^{-1}|Q|$ for any one of these squares. □

## 8. Consequences of Lemma 7.2: Planiness and graininess

In this section we shall utilize the first conclusion Lemma 7.2 to enforce a "planiness" property on the tubes in $\mathbb{T}_\rho$; roughly speaking, we require for all $x_0$ satisfying the property $P_2$ that most of the tubes in $\mathbb{T}_\rho(x_0)$ to be coplanar. Then we will shift scales in order to also enforce this property on the $\delta$-tubes $\mathbb{T}$. Later in this section we will also use the second conclusion of Lemma 7.2 to also enforce a certain "graininess" property on the set $E_\delta$.

We begin this program with

PROPOSITION 8.1. *Let $x_0$ be such that $P_2(x_0)$ holds. Then one can find a family of planes $\Pi_\rho = \Pi_\rho(x_0)$ through $x_0$ and subcollections $\mathbb{T}_\rho(x_0, \pi) \subset \mathbb{T}_\rho(x_0)$ for all $\pi \in \Pi_\rho(x_0)$ such that*

- *For all $\pi \in \Pi_\rho(x_0)$ and $T_\rho \in \mathbb{T}_\rho(x_0, \pi)$ we have*

$$T_\rho \subset N(\pi), \tag{66}$$

*where $N(\pi)$ is the $\delta^{-C_2\sqrt{\epsilon}}\rho$-neighbourhood of $\pi$, and $C_2$ is a constant to be chosen later.*

- *The $\mathbb{T}_\rho(x_0, \pi)$ are disjoint as $\pi$ varies, and*

$$\#\mathbb{T}_\rho(x_0, \pi) \gtrapprox \delta^{C\sqrt{\epsilon}}\rho^{-(3-d)} \tag{67}$$

*for all $\pi \in \Pi_\rho(x_0)$.*

- *We have*

$$1 \leq \#\Pi_\rho(x_0) \lessapprox \delta^{-C\sqrt{\epsilon}}. \tag{68}$$

- *We have*

$$\mathop{\tilde{\forall}}_{T_\rho} T_\rho \in \mathbb{T}_\rho(x_0) : T_\rho \in \bigcup_{\pi \in \Pi_\rho(x_0)} \mathbb{T}_\rho(x_0, \pi) \tag{69}$$

- *For every $\pi \in \Pi_\rho(x_0)$, the set $\mathbb{T}_\rho(x_0, \pi)$ is a subset of*

(70) $$\{T_\rho \in \mathbb{T}_\rho(x_0) : P_1(x_0, T_\rho)\}.$$

For the purposes of visualization we recommend thinking of $\Pi_\rho$ as consisting of a single plane $\pi$, and $\mathbb{T}_\rho(x_0, \pi)$ as equal to $\mathbb{T}_\rho(x_0)$.

*Proof.* The idea of the proof is to iterate the following lemma.

LEMMA 8.2. *Let $x_0$ be such that $P_2(x_0)$ holds. Suppose that $\mathbb{T}'_\rho$ is a subset of* (70) *such that*

(71) $$\#\mathbb{T}'_\rho \geq \delta^{\sqrt{\epsilon}} \#\mathbb{T}_\rho(x_0).$$

*Then one can find a plane $\pi$ through $x_0$ such that*

(72) $$\#\{T_\rho \in \mathbb{T}'_\rho : T_\rho \subset N(\pi)\} \gtrapprox \delta^{C_2\sqrt{\epsilon}} \rho^{-(3-d)},$$

*if $C_2$ is chosen sufficiently large.*

*Proof.* Suppose for contradiction that (72) failed for all planes $\pi$ through $x_0$. Let $F$ denote the set

(73) $$F = \{x \in E_\delta \cap B(x_0, C\rho) : (42) \text{ holds for all } \delta \leq \sigma \ll 1\}.$$

Let $T_\rho$ be any element of $\mathbb{T}'_\rho$. From (45), (46), and Lemma 5.2 we have

(74) $$\tilde{\forall}_{T,x} T \in \mathbb{T}[T_\rho], x \in T \cap B(x_0, C\rho) : x \in F, \mu_\delta[T_\rho](x) \leq \delta^{-\sqrt{\epsilon}} \rho^{-(3-d)}.$$

From (31) we thus obtain

$$\int_{F \cap \{\mu_\delta[T_\rho] \leq \delta^{-\sqrt{\epsilon}} \rho^{-(3-d)}\}} \mu_\delta[T_\rho] \approx \int_{B(x_0, C\rho)} \mu_\delta[T_\rho].$$

By (27) and (11) we see that

$$\delta^{-\sqrt{\epsilon}} \rho^{-(3-d)} |A(x_0, T_\rho) \cap F| \gtrapprox \delta^{3/2},$$

where $A(x_0, T_\rho)$ was defined in (57). We rearrange this as

$$\int_F \chi_{A(x_0, T_\rho)} \gtrapprox \delta^{\sqrt{\epsilon}} \delta^{3-d} \rho^d.$$

Summing this over all $T_\rho \in \mathbb{T}'_\rho$ and using (71) we obtain

(75) $$\int_F \sum_{T_\rho \in \mathbb{T}'_\rho} \chi_{A(x_0, T_\rho)} \gtrapprox \delta^{C\sqrt{\epsilon}} \#\mathbb{T}_\rho(x_0)[\delta^{3-d} \rho^d].$$

On the other hand, we will shortly prove the estimate

(76)
$$\sum_{T_\rho^1}\sum_{T_\rho^2}\sum_{T_\rho^3} |A(x_0, T_\rho^1) \cap A(x_0, T_\rho^2) \cap A(x_0, T_\rho^3) \cap F| \lessapprox \delta^{-C\sqrt{\epsilon}} \delta^{cC_2\sqrt{\epsilon}} \delta^{3-d} \rho^d \#\mathbb{T}_\rho(x_0)^3,$$

where the tubes $T_\rho^i$ will always be assumed to range over $\mathbb{T}_\rho'$, and $c$ denotes a constant independent of $C_2$ and $\varepsilon$ but which can vary from line to line. Rewriting the left-hand side of (76) as

$$\Big\| \sum_{T_\rho \in \mathbb{T}_\rho'} \chi_{A(x_0, T_\rho)} \Big\|_{L^3(F)}^3$$

and using (8) (with $d$ replaced by 3/2) and (51), we thus obtain

$$\Big\| \sum_{T_\rho \in \mathbb{T}_\rho'} \chi_{A(x_0, T_\rho)} \Big\|_{L^1(F)} \lessapprox \delta^{-C\sqrt{\epsilon}} \delta^{cC_2\sqrt{\epsilon}} \delta^{3-d} \rho^d \#\mathbb{T}_\rho(x_0)$$

which contradicts (75) if $C_2$ is chosen sufficiently large.

It remains to prove (76). Fix $T_\rho^1, T_\rho^2, T_\rho^3$, and choose $v_i$ so that (59) holds. We divide the contributions to (76) into three cases.

First, we consider the contribution when $v_1$ and $v_2$ are close together, or more precisely when

(77)
$$\angle v_1, v_2 \lessapprox \delta^{C_2\sqrt{\epsilon}/2}.$$

Estimating the summand crudely by (51), we may majorize the contribution of this case by

(78) $\qquad \delta^{-C\sqrt{\epsilon}} \delta^{3-d} \rho^d \#\{(T_\rho^1, T_\rho^2, T_\rho^3) \in {\mathbb{T}_\rho'}^3 : (77) \text{ holds}\}.$

However, by (52) we see that for any fixed $T_\rho^1$ there are at most

$$\lessapprox \delta^{-C\sqrt{\epsilon}} (\delta^{C_2\sqrt{\epsilon}/2})^c \rho^{-(3-d)}$$

tubes $T_\rho^2$ for which (77) holds; by (49) this is majorized by

$$\lessapprox \delta^{-C\sqrt{\epsilon}} \delta^{cC_2\sqrt{\epsilon}} \#\mathbb{T}_\rho(x_0).$$

Thus the contribution of this case to (76) is acceptable.

Next we estimate the contribution for which (77) fails, but for which

(79)
$$\angle \pi(v_1, v_2), v_3 \lessapprox \delta^{-C_2\sqrt{\epsilon}} \rho.$$

where $\pi(v_1, v_2)$ is the plane spanned by $v_1, v_2$. By crudely estimating the summand by (51) again, we see that the contribution to (76) is at most

$$\delta^{-C\sqrt{\epsilon}} \delta^{3-d} \rho^d \#\{(T_\rho^1, T_\rho^2, T_\rho^3) \in (\mathbb{T}_\rho')^3 : (79) \text{ holds}\}.$$

However, we are assuming (72) to fail for all planes $\pi$. This implies that for any fixed $T_\rho^1$, $T_\rho^2$, there are at most $\delta^{C_2\sqrt{\epsilon}} \rho^{-(3-d)}$ tubes $T_\rho^3$ for which (79) holds. By repeating the previous calculation we see that the contribution of this case to (76) is also acceptable.

Finally, we consider the contribution for which (77) and (79) both fail. In this case the first part Lemma 7.3 applies with $M = \delta^{-C_2\sqrt{\epsilon}/2}$. Using (60) we therefore see that this contribution is also acceptable. Combining the three contributions together we obtain the desired estimate. $\square$

With this lemma it is easy to prove the proposition by the following construction. Initialize $\Pi_\rho(x_0)$ to be the empty set and $\mathbb{T}'_\rho$ to be (70). From (50) we see that (71) currently holds.

We now perform the following iteration. We use Lemma 8.2 to extract a plane $\pi$ which obeys (72) with respect to $\mathbb{T}'_\rho$. We add this plane $\pi$ to $\Pi_\rho(x_0)$ and define
$$\mathbb{T}_\rho(x_0, \pi) = \{T_\rho \in \mathbb{T}'_\rho : T_\rho \subset N(\pi)\}.$$
We now remove all the elements of $\mathbb{T}_\rho(x_0, \pi)$ from $\mathbb{T}'_\rho$ and repeat the iteration until (71) fails, at which point we halt the construction.

It is easy to verify that this iteration gives all the properties claimed in the proposition. The bound (69) follows by construction. The upper bound on $\#\Pi_\rho$ follows from (72) and (49), while the lower bound on $\#\Pi_\rho$ follows from (69). $\square$

Henceforth we assume that the collections $\Pi_\rho(x_0)$ and $\mathbb{T}_\rho(x_0, \pi)$ have been chosen for all $x_0$ such that $P_2(x_0)$ held.

Having obtained planiness, we now turn to the companion property of graininess. We remark that the heuristic that one should be able to reduce the study of Besicovitch sets to the grainy case is due to Tom Wolff (private communication). Much of the work in this paper was motivated by a desire to make this heuristic rigorous.

*Definition* 8.3. Let $\pi$ be a plane, and $x_0$ be a point in $\mathbb{R}^3$. A set $S$ is said to be *grainy* with respect to $(x_0, \pi)$ if $S \subset B(x_0, C\rho)$ and
$$(80) \qquad S = \bigcup_{Q \in \mathbb{Q}} Q$$
for some collection $\mathbb{Q}$ of squares $Q$ parallel to $\pi$ such that
$$(81) \qquad \#\mathbb{Q} \lessapprox \delta^{-C\sqrt{\epsilon}} \delta^{(2-d)/2}.$$
If $S$ is a grainy set, we define its dilate $CS$ to be
$$CS = \bigcup_{Q \in \mathbb{Q}} CQ$$
where $CQ$ is the dilate of $Q$ by $C$ around its center.

Note if $S$ is grainy, then $CS$ is also grainy with slightly larger constants. Also, from (81) we have

$$(82) \qquad |S|, |CS| \lesssim \delta\rho^2 \#\mathbb{Q} \lessapprox \delta^{-C\sqrt{\epsilon}} \delta^{3-d} \rho^d.$$

PROPOSITION 8.4. *Let the hypotheses and notation be as in Proposition* 8.1. *Then for each* $\pi \in \Pi_\rho(x_0)$ *there exists a set* $S(x_0, \pi)$ *which is grainy with respect to* $(x_0, \pi)$, *such that*

$$(83) \qquad \tilde{\forall}_{T_\rho, T, x} T_\rho \in \mathbb{T}_\rho(x_0, \pi), T \in \mathbb{T}[T_\rho], x \in T \cap B(x_0, C\rho) : x \in S(x_0, \pi).$$

*Proof.* We will modify the argument in Proposition 8.1.

Fix $x_0$, $\pi$. Cover $B(x_0, C\rho)$ by a finitely overlapping collection of squares $Q$ parallel to $\pi$. Define $\mathbb{Q}(x_0, \pi)$ to be the subcollection of these squares $Q$ for which

$$(84) \qquad \frac{|E_\delta \cap Q|}{|Q|} \gtrapprox \delta^{C_3 \sqrt{\epsilon}}.$$

where $C_3 > 0$ is a constant we will choose later. Define $S(x_0, \pi)$ by

$$S(x_0, \pi) = \bigcup_{Q \in \mathbb{Q}(x_0, \pi)} Q.$$

From (84), (51) and the estimate $|Q| \sim \delta^2$ we obtain (81), and so $S(x_0, \pi)$ is grainy with respect to $(x_0, \pi)$.

Define the set $X$ by

$$X = (E_\delta \cap B(x_0, C\rho)) \setminus S(x_0, \pi).$$

Then by construction

$$(85) \qquad \sup_Q \frac{|X \cap Q|}{|Q|} \lessapprox \delta^{C_3\sqrt{\epsilon}}$$

where $Q$ ranges over all squares parallel to $\pi$.

To finish the proof we must show (83). By summing (46) for all $T_\rho \in \mathbb{T}_\rho(x_0, \pi)$ we have

$$\tilde{\forall}_{T_\rho, T, x} T_\rho \in \mathbb{T}_\rho(x_0, \pi), T \in \mathbb{T}[T_\rho], x \in T \cap B(x_0, C\rho) : (47) \text{ holds}.$$

Thus by Lemma 5.2 it suffices to show that

$$\tilde{\forall}_{T_\rho} T_\rho \in \mathbb{T}_\rho(x_0, \pi), T \in \mathbb{T}[T_\rho], x \in T \cap B(x_0, C\rho) : ((47) \Longrightarrow x \notin X).$$

We rewrite this as

$$(86) \qquad \int_X \sum_{T_\rho \in \mathbb{T}_\rho(x_0, \pi) : \mu_\delta[T_\rho] \leq \delta^{-\sqrt{\epsilon}} \rho^{-(3-d)}} \mu_\delta[T_\rho] \lessapprox \delta^{c\sqrt{\epsilon}} \int_{B(x_0, C\rho)} \sum_{T_\rho \in \mathbb{T}_\rho(x_0, \pi)} \mu_\delta[T_\rho].$$

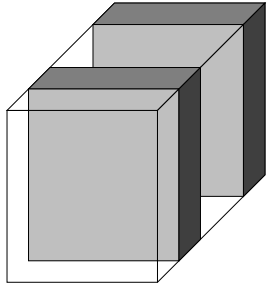

Figure 5. A (simplified) typical intersection of $E_\delta$ with a $\rho$-cube. In this case the plane $\pi$ is oriented in the xz direction.

We majorize the left-hand side of (86) as

$$\int_X \sum_{T_\rho \in \mathbb{T}_\rho(x_0,\pi)} \delta^{-\sqrt{\epsilon}} \rho^{-(3-d)} \chi_{A(x_0,T_\rho)}.$$

The right-hand side of (86) can be estimated using (67) and (27). Thus (86) reduces to showing that

$$\int_X \sum_{T_\rho \in \mathbb{T}_\rho(x_0,\pi)} \delta^{-\sqrt{\epsilon}} \rho^{-(3-d)} \chi_{A(x_0,T_\rho)} \lessapprox \delta^{C\sqrt{\epsilon}} \rho^{-(3-d)} \delta^{3/2}$$

for some sufficiently large $C$. This simplifies by (49) to

$$\int_X \sum_{T_\rho \in \mathbb{T}_\rho(x_0,\pi)} \chi_{A(x_0,T_\rho)} \lessapprox \delta^{C\sqrt{\epsilon}} \delta^{3-d} \rho^d \# \mathbb{T}_\rho(x_0).$$

We will shortly prove

(87) $$\sum_{T_\rho^1} \sum_{T_\rho^2} \sum_{T_\rho^3} |A(x_0, T_\rho^1) \cap A(x_0, T_\rho^2) \cap A(x_0, T_\rho^3) \cap X|$$
$$\lessapprox \delta^{-C\sqrt{\epsilon}} \delta^{cC_3\sqrt{\epsilon}} \delta^{3-d} \rho^d \# \mathbb{T}_\rho(x_0)^3,$$

where the $T_\rho^i$ are assumed to range over $\mathbb{T}_\rho(x_0, \pi)$; the desired estimate then follows as before from (8) with $d$ replaced by $3/2$ and then (51), if $C_3$ is chosen sufficiently large.

As in the proof of Proposition 8.1, the contributions to (87) when

(88) $$\angle v_1, v_2 \lessapprox \delta^{C_3\sqrt{\epsilon}}$$

or when

(89) $$\angle \pi(v_1, v_2), v_3 \gtrapprox \delta^{-2C_3\sqrt{\epsilon}} \rho.$$

are acceptable if $C_3$ is sufficiently large. Thus it only remains to consider the contribution for which (88) and (89) both fail. However, in this case the second part of Lemma 7.2 applies. From (61) and (85) we have

$$|A(x_0, T_\rho^1) \cap A(x_0, T_\rho^2) \cap A(x_0, T_\rho^3) \cap X| \lessapprox \delta^{-C\sqrt{\epsilon}} \delta^{3-d} \rho^d \delta^{cC_3\sqrt{\epsilon}}.$$

Thus the contribution from this case to (87) are also acceptable, as desired. □

Henceforth we assume that the sets $S(x_0, \pi)$ have been chosen for all $x_0$ such that $P_2(x_0)$ holds, and all $\pi \in \Pi_\rho(x_0)$. The constants $C_2$ and $C_3$ are now fixed, and future constants $C$ may implicitly depend on these two.

## 9. Consistency of planes and squares

In Proposition 6.6 we showed that any sticky collection of tubes satisfying (20) and (21) had to satisfy the property $P_2(x)$ extremely often. This property in turn implied two strong properties on the collection, namely planiness and graininess, by Propositions 8.1 and 8.4.

The planiness property is really a property about the $\rho$-tubes $\mathbb{T}_\rho$ rather than the $\delta$-tubes $\mathbb{T}$. By taking advantage of the double-sticky assumption and changing scale from $\delta$ to $\delta^2$, we shall be able to also impose a planiness property in the $\delta$-tubes. More precisely, we show

*Definition* 9.1. Define $P_3(x)$ to be property that one can find a family of planes $\Pi_\delta(x)$ through $x$ and subcollections $\mathbb{T}(x, \pi) \subset \mathbb{T}(x)$ for all $\pi \in \Pi_\delta(x)$ such that:

- For all $\pi \in \Pi_\delta(x)$ and $T \in \mathbb{T}(x, \pi)$ we have

(90) $$T \subset N_{\delta^{-C\sqrt{\epsilon}}\delta}(\pi).$$

- The $\mathbb{T}(x, \pi)$ are disjoint as $\pi$ varies, and $\#\mathbb{T}(x, \pi) \gtrapprox \delta^{C\sqrt{\epsilon}} \delta^{d-3}$ for all $\pi \in \Pi_\delta(x)$.

- We have

(91) $$1 \leq \#\Pi_\delta(x) \lessapprox \delta^{-C\sqrt{\epsilon}}.$$

- We have

(92) $$\tilde{\forall}_T T \in \mathbb{T}(x) : T \in \mathbb{T}(x, \pi) \text{ for some } \pi \in \Pi_\delta(x).$$

- We have

(93) $$\#\{T \in \mathbb{T}(x_0) : \mathrm{dir}(T) \in B(\omega, \theta)\} \lessapprox \delta^{-C\sqrt{\epsilon}} \theta^c \delta^{d-3}$$

for all directions $\omega$ and all $\delta \leq \theta \ll 1$.

PROPOSITION 9.2. *Suppose that there is an X-ray estimate at dimension $d$, and that there exists a Besicovitch set $E$ with $\overline{\dim}(E) < d + \varepsilon$. Then there exists a sticky collection $\mathbb{T}$ of tubes such that (20) and (21) hold, and such that*

(94) $$\widetilde{\forall}_{T,x} T \in \mathbb{T}, x \in T : P_3(x).$$

*Proof.* By Proposition 3.3 with $\delta$ replaced by $\delta^2$, we can find a doubly sticky collection $\mathbb{T}_{\delta^2}$ of tubes with parent collection $\mathbb{T}_\delta$ and grandparent collection $\mathbb{T}_\rho$ which obey the analogues of (20), (21), and (22) with $\delta$ replaced by $\delta^2$. By Proposition 6.6 the property $P_2$ (with $\delta$ replaced by $\delta^2$) holds almost always, and the conclusions of Proposition 8.1 obtain whenever $P_2$ is satisfied. Similarly the property (52) with $\delta$ replaced by $\delta^2$ also holds almost always. By comparing these statements against Definition 9.1 we see that the desired properties obtain with $\mathbb{T} = \mathbb{T}_\delta$. □

Henceforth we fix $\mathbb{T}$ to be a collection of tubes satisfying all the conclusions of Proposition 9.2, and assume that we have chosen the collections $\Pi_\delta(x)$ and $\mathbb{T}(x, \pi)$ for all $x$ such that $P_2(x)$ holds, and all $\pi \in \Pi_\delta(x)$. This changing of scale from $\delta$ to $\rho$ is intriguingly reminiscent of the techniques employed in proving restriction theorems (see e.g. [1], [2], [3], [16], [14]).

Since $\mathbb{T}$ is sticky and satisfies (20) and (21), we may apply Proposition 6.6 again and obtain (53). In particular, the conclusions of Propositions 8.1 and 8.4 obtain for a large number of $x_0$.

We thus have two scales of structure on our set, one at scale $\delta$ arising from $P_3$, and one arising from the propositions in the previous section. Fortunately, there is a close relationship between the two structures, namely that the planes in $\Pi_\delta(x)$ and the planes in $\Pi_\rho(x_0)$ are mostly parallel if $x \in B(x_0, C\rho)$. In particular, the planes in $\Pi_\delta(x)$ are mostly parallel to the squares $Q$ that contain $x$. More precisely:

*Definition* 9.3. We define $P_4(x_0, x, T_\rho, T)$ to be the property that $P_2(x_0)$ and $P_3(x)$ hold, that $x \in B(x_0, C\rho)$, and there exist planes $\pi \in \Pi_\delta(x)$, $\pi' \in \Pi_\rho(x_0)$ such that $T \in \mathbb{T}(x, \pi) \cap \mathbb{T}[T_\rho]$, $T_\rho \in \mathbb{T}_\rho(x_0, \pi')$, $x \in S(x_0, \pi')$, and we have the consistency condition

(95) $$\angle \pi, \pi' \lessapprox \delta^{-C\sqrt{\epsilon}}\rho.$$

PROPOSITION 9.4. *If the constants in the above definition are chosen appropriately, then we have*
(96)
$$\widetilde{\forall}_{T_\rho, x_0, T, x} T_\rho \in \mathbb{T}_\rho, x_0 \in T_\rho, T \in \mathbb{T}[T_\rho], x \in T \cap B(x_0, C\rho) : P_4(x_0, x, T_\rho, T).$$

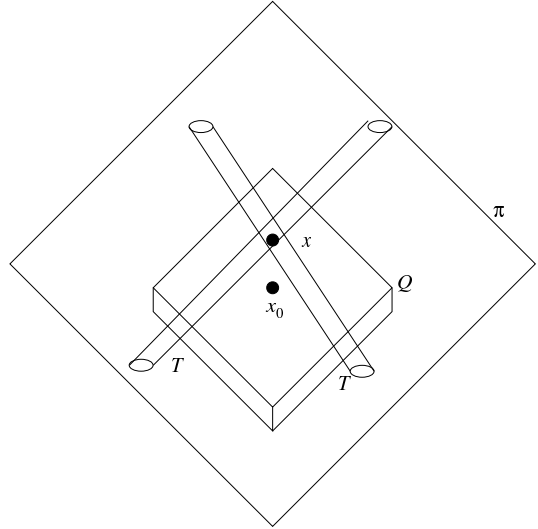

Figure 6. A typical situation in which $P_4(x_0, x, T_\rho, T)$ holds. The tubes through $x$ lie on the plane $\pi$, which is essentially parallel to the square $Q \subset S(x_0, \pi')$ which contains $x$. The plane $\pi'$ (not pictured) differs in angle from $\pi$ by about $\rho$.

Unfortunately, the proof of this proposition is the most technical part of the paper, and mostly consists of the combinatorial reshuffling machinery developed in Section 5. Heuristically, the idea is as follows. If the planes in $\Pi_\delta(x)$ are at a large angle to the planes in $\Pi_\rho(x_0)$, then any tubes $T$ which are parallel to both tubes will be constrained to a narrow angular band. By (52) this would imply that very few tubes are parallel to both kinds of planes, which will eventually contradict (66), (50) and (90), (69).

*Proof.* We need to show that $\tilde{\forall}_{T_\rho, x_0, T, x} Q : P_4$, where $Q = Q(x_0, x, T_\rho, T)$ denotes the property that

$$T_\rho \in \mathbb{T}_\rho, x_0 \in T_\rho, T \in \mathbb{T}[T_\rho], x \in T \cap B(x_0, C\rho).$$

Apart from the condition (95), this can be achieved simply by applying the machinery of Section 5 to the previous propositions, as follows.

From (69) and Lemma 5.5 we have

$$\tilde{\forall}_{x_0, T_\rho} P_2(x_0), T_\rho \in \mathbb{T}_\rho(x_0) : T_\rho \in \bigcup_{\pi' \in \Pi_\rho(x_0)} \mathbb{T}_\rho(x_0, \pi').$$

Clearly this implies

$$\tilde{\forall}_{x_0, T_\rho} P_2(x_0), T_\rho \in \mathbb{T}_\rho(x_0) : \Big[\tilde{\forall}_{T, x} T \in \mathbb{T}[T_\rho], x \in T : T_\rho \in \bigcup_{\pi' \in \Pi_\rho(x_0)} \mathbb{T}_\rho(x_0, \pi')\Big].$$

By Lemma 5.10 we therefore see that

(97) $$\tilde{\forall}_{x_0, x, T_\rho, T} P_2(x_0), Q : T_\rho \in \bigcup_{\pi' \in \Pi_\rho(x_0)} \mathbb{T}_\rho(x_0, \pi').$$

On the other hand, from (53) and another application of Lemma 5.6 we have

$$\tilde{\forall}_{x_0, x, T_\rho, T} Q : P_2(x_0).$$

Combining this with (97) using Lemma 5.4 we obtain

(98) $$\tilde{\forall}_{x_0, x, T_\rho, T} Q : P_2(x_0), T_\rho \in \bigcup_{\pi' \in \Pi_\rho(x_0)} \mathbb{T}_\rho(x_0, \pi').$$

Now from (83) and Lemma 5.5 we have

$$\tilde{\forall}_{x_0, \pi', T_\rho, T, x} P_2(x_0), \pi' \in \Pi_\rho(x_0), T_\rho \in \mathbb{T}_\rho(x_0, \pi'),$$
$$T \in \mathbb{T}[T_\rho], x \in T \cap B(x_0, C\rho) : x \in S(x_0, \pi'),$$

which we may rearrange using Corollary 5.7 as

$$\tilde{\forall}_{x_0, x, T_\rho, T} Q, P_2(x_0), T_\rho \in \bigcup_{\pi' \in \Pi_\rho(x_0)} \mathbb{T}_\rho(x_0, \pi') : x \in S(x_0, \pi')$$

with the understanding that $\pi'$ is the unique plane in $\Pi_\rho(x_0)$ such that $T_\rho \in \mathbb{T}_\rho(x_0, \pi')$. Combining this with (98) via Lemma 5.4, we obtain

(99) $$\tilde{\forall}_{x_0, x, T_\rho, T} Q : R$$

where $R = R(x_0, x, T_\rho, T)$ denotes the property that $P_2(x_0)$ holds, and for some $\pi' \in \Pi_\rho(x_0)$ we have $T_\rho \in \mathbb{T}_\rho(x_0, \pi')$ and $x \in S(x_0, \pi')$.

From (92) and Lemma 5.5 we have

$$\tilde{\forall}_{x, T} P_3(x), T \in \mathbb{T}(x) : T \in \bigcup_{\pi \in \Pi_\delta(x)} \mathbb{T}(x, \pi).$$

We rearrange this using Corollary 5.7 as

$$\tilde{\forall}_{T_\rho, T, x} T_\rho \in \mathbb{T}_\rho, T \in \mathbb{T}[T_\rho], x \in T, P_3(x) : T \in \bigcup_{\pi \in \Pi_\delta(x)} \mathbb{T}(x, \pi).$$

By (94) and Lemma 5.4 we thus have

$$\tilde{\forall}_{T_\rho, T, x} T_\rho \in \mathbb{T}_\rho, T \in \mathbb{T}[T_\rho], x \in T : R',$$

where $R' = R'(x, T)$ denotes the property that

$$P_3(x), \text{ and } T \in \bigcup_{\pi \in \Pi_\delta(x)} \mathbb{T}(x, \pi).$$

By Lemma 5.10 we then obtain

(100) $$\tilde{\forall}_{x_0, x, T_\rho, T} Q : R'.$$

The property $P_4(x_0, x, T_\rho, T)$ is equivalent to

$$R, R', (95)$$

where $C_4$ is a constant to be chosen later, and we adopt the convention that $\pi$ is the unique plane in $\Pi_\delta(x)$ such that $T \in \mathbb{T}(x, \pi)$. Thus by (99), (100), and Corollary 5.3 it only remains to show that

$$\tilde{\forall}_{x_0, x, T_\rho, T} Q : ((R, R') \Longrightarrow (95)),$$

or in other words that
(101)
$$|\{(x_0, x, T_\rho, T) : Q, R, R', \angle \pi, \pi' \gtrapprox \delta^{-C_4\sqrt{\epsilon}}\rho\}| \lessapprox \delta^{c\sqrt{\epsilon}}|\{(x_0, x, T_\rho, T) : Q\}|.$$

The left-hand side of (101) is majorized by

(102) $$|\{(x_0, x, T_\rho, T, \pi, \pi') : P_2(x_0), P_3(x), \pi' \in \Pi_\rho(x_0), \pi \in \Pi_\delta(x), |x - x_0| \leq \rho,$$
$$T_\rho \in \mathbb{T}_\rho(x_0, \pi'), T \in \mathbb{T}[T_\rho] \cap \mathbb{T}(x, \pi), \angle \pi, \pi' \gtrapprox \delta^{-C_4\sqrt{\epsilon}}\rho\}|.$$

Of course, $T_\rho, T, \pi, \pi'$ are measured by counting measure and $x_0, x$ by Lebesgue measure.

Fix $x_0, x, \pi, \pi'$ such that

$$P_2(x_0), P_3(x), \pi' \in \Pi_\rho(x_0), \pi \in \Pi_\delta(x), \angle \pi, \pi' \gtrapprox \delta^{-C_4\sqrt{\epsilon}}\rho,$$

and suppose that $T_\rho$ and $T$ be such that $T_\rho \in \mathbb{T}_\rho(x_0, \pi')$ and $T \in \mathbb{T}[T_\rho] \cap \mathbb{T}(x, \pi)$. From (66) and (90) we see that $T$ is contained in the $\delta^{-C\sqrt{\epsilon}}\rho$-neighbourhood of $\pi'$, and the $\delta^{-C\sqrt{\epsilon}}\delta$-neighbourhood of $\pi$. From the angular separation assumption on $\pi$ and $\pi'$ and elementary geometry this implies that the direction of $T$ is constrained to a cap of diameter $O(\delta^{-C\sqrt{\epsilon}}\delta^{C_4\sqrt{\epsilon}})$, whose center depends only on $\pi$ and $\pi'$. From (93) we thus see that the contribution of $(x_0, x, \pi, \pi')$ to (102) is majorized by

$$\delta^{-C\sqrt{\epsilon}}(\delta^{-C\sqrt{\epsilon}}\delta^{C_4\sqrt{\epsilon}})^c\delta^{d-3}.$$

We may therefore majorize (102) by

$$\int\int_{P_2(x_0), P_3(x), |x-x_0| \leq C\rho} \sum_{\pi' \in \Pi_\rho(x_0)} \sum_{\pi \in \Pi_\delta(x)} \delta^{-C\sqrt{\epsilon}}(\delta^{-C\sqrt{\epsilon}}\delta^{C_4\sqrt{\epsilon}})^c\delta^{d-3} \, dx dx_0.$$

From (68) and (91) this is majorized by

$$\delta^{-C\sqrt{\epsilon}}\delta^{cC_4\sqrt{\epsilon}}\delta^{d-3} \int_{P_2(x_0)} \left( \int_{E_\delta \cap B(x_0, C\rho)} dx \right) dx_0.$$

By (51), the inner integral is $\lessapprox \delta^{-C\sqrt{\epsilon}}\delta^{3-d}\rho^d$. From (21) we can therefore majorize this as

(103) $$\delta^{-C\sqrt{\epsilon}}\delta^{cC_4\sqrt{\epsilon}}\delta^{3/2}.$$

We now turn our attention to the right-hand side of (101), which we can rewrite as

$$\delta^{c\sqrt{\epsilon}} \int\int_{|x-x_0|\leq C\rho} \sum_{T_\rho\in\mathbb{T}_\rho} \sum_{T\in\mathbb{T}[T_\rho]} \chi_{T_\rho}(x_0)\chi_T(x)\,dxdx_0$$

for some $c > 0$. Evaluating the $x$-integral using (27), we see that this is

$$\approx \delta^{c\sqrt{\epsilon}} \int \sum_{T_\rho\in\mathbb{T}_\rho} \chi_{T_\rho}(x_0)\delta^{3/2}\,dx_0,$$

which by (25) is

$$\approx \delta^{c\sqrt{\epsilon}}\delta^{3/2}.$$

Comparing this with (103) we obtain (101) as desired, if $C_4$ is chosen sufficiently large. $\square$

Applying Lemma 5.10 to this proposition, we obtain the main result of this section, which is phrased in terms of a rather complicated property $P_5$.

COROLLARY 9.5. *If $P_5(x_0, T_\rho)$ denotes the property that*

(104) $$\tilde{\forall}_{T,x} T \in \mathbb{T}[T_\rho], x \in T \cap B(x_0, C\rho) : P_4(x_0, x, T),$$

*then one has*

(105) $$\tilde{\forall}_{T_\rho, x_0} T_\rho \in \mathbb{T}_\rho, x_0 \in T_\rho : P_5(x_0, T_\rho).$$

## 10. An arithmetic progression of length three

In the last few sections we have been steadily accumulating structural information on a collection $\mathbb{T}$ of tubes. This program has culminated in Corollary 9.5, in which we have shown that a rather complicated property $P_5$ holds almost always.

We have taken great pains to ensure that this property $P_5$ almost never failed. We will now see the payoff for this carefulness, in that we will be easily able to find three widely spaced points in arithmetic progression for which $P_5$ holds. Furthermore we shall be able to impose an additional technical condition, namely that the planes associated to these points are not too parallel. More precisely:

PROPOSITION 10.1. *Let $\mathbb{T}$ be a sticky collection of tubes satisfying (105). Then there exists a tube $T_\rho \in \mathbb{T}_\rho$ and three points $x_0, x_0 + s, x_0 + 2s$ in $T_\rho$ such that $P_5(x_0 + si, T_\rho)$ holds for $i = 0, 1, 2$. Furthermore, $|s| \sim 1$ and $\angle \pi(x_0 + si), \pi(x_0 + sj) \gtrapprox \delta^{C_5 \sqrt{\epsilon}}$ for distinct $i, j \in \{0, 1, 2\}$, where $\pi(x)$ denotes the plane for which $T_\rho \in \mathbb{T}_\rho(x, \pi(x))$, and $C_5 > 0$ is a suitable constant.*

Note that the plane $\pi(x)$ is well defined if $P_5(x, T_\rho)$ holds.

*Proof.* From (105) and Lemma 5.9 we have

$$\tilde{\forall}_{T_\rho} T_\rho \in \mathbb{T}_\rho : [\tilde{\forall}_x x \in T_\rho : P_5(x, T_\rho)].$$

Thus by (31) we can find a tube $T_\rho \in \mathbb{T}_\rho$ such that

(106) $$\tilde{\forall}_x x \in T_\rho : P_5(x, T_\rho).$$

Fix such a tube $T_\rho$, and let $A$ denote the set

$$A = \left\{ (x_0, s) \in \mathbb{R}^3 \times \mathbb{R}^3 : |s| \sim 1, \quad x_0 + si \in T_\rho \text{ for } i = 0, 1, 2 \right\}.$$

Note that $x_0$ lives in $T_\rho$, and $s$ lives in a translate of $CT_\rho$. We have to show that the set

$$\left\{ (x_0, s) \in A : P_5(x_0 + si, T_\rho), \angle \pi(x_0 + si), \pi(x_0 + sj) \gtrapprox \delta^{C_5 \sqrt{\epsilon}} \right.$$
$$\left. \text{for all distinct } i, j = 0, 1, 2 \right\}$$

is nonempty. On the other hand, from (106) we have

$$|\{(x_0, s) \in A : P_5(x_0 + si, T_\rho) \text{ fails}\}| \lessapprox \delta |\{x \in T_\rho : P_5(x, T_\rho) \text{ fails}\}| \lessapprox \delta^{c\sqrt{\epsilon}} \delta^2$$

for all $i = 0, 1, 2$. Since $|A| \sim \delta^2$, it will therefore suffice to show that

$$|\{(x_0, s) \in A : P_5(x_0 + si, T_\rho), P_5(x_0 + sj, T_\rho), \angle \pi(x_0 + si), \pi(x_0 + sj) \lessapprox \delta^{C_5 \sqrt{\epsilon}}\}| \lessapprox \delta^{c\sqrt{\epsilon}} \delta^2$$

for all distinct $i, j = 0, 1, 2$. Making the change of variables $x = x_0 + si$ and $y = x_0 + sj$, it therefore suffices to show that

(107) $$|\{(x, y) \in F \times F : \angle(\pi(x), \pi(y)) \lessapprox \delta^{C_5 \sqrt{\epsilon}}\}| \lessapprox \delta^{c\sqrt{\epsilon}} \delta^2$$

for some $c > 0$, where $F$ denotes the set

$$F = \{x \in T_\rho : P_5(x, T_\rho)\}.$$

The idea is as follows. By using the Kakeya estimate (14) one can easily recover the trivial estimate for (107) with $c = 0$. To improve upon this we use the X-ray estimate (15) instead, observing that if the squares corresponding to different $x$ are almost parallel then there will be many parallel tubes in each direction.

We now turn to the details. Let $\Pi$ be a maximal $\delta^{C_5\sqrt{\epsilon}}$-angular-separated set of planes containing the axis of $T_\rho$. If $x \in F$, then there exists a $\pi \in \Pi$ such that $\angle(\pi(x), \pi) \lessapprox \delta^{C_5\sqrt{\epsilon}}$. We can therefore majorize the left-hand side of (107) by

$$(107) \lessapprox \sum_{\pi \in \Pi} \left|\left\{(x,y) \in F \times F : \angle(\pi(x), \pi), \angle(\pi(y), \pi) \lessapprox \delta^{C_5\sqrt{\epsilon}}\right\}\right|$$

which we rewrite as

$$(108) \qquad (107) \lessapprox \sum_{\pi \in \Pi} \left|\left\{x \in F : \angle(\pi(x), \pi) \lessapprox \delta^{C_5\sqrt{\epsilon}}\right\}\right|^2.$$

On the other hand, we clearly have

$$\sum_{\pi \in \Pi} \left|\left\{x \in F : \angle(\pi(x), \pi) \lessapprox \delta^{C_5\sqrt{\epsilon}}\right\}\right| \lessapprox |F| \lessapprox \delta.$$

It thus suffices to show

$$(109) \qquad \left|\left\{x \in F : \angle(\pi(x), \pi) \lessapprox \delta^{C_5\sqrt{\epsilon}}\right\}\right| \lessapprox \delta^{cC_5\sqrt{\epsilon}} \delta^{-C\sqrt{\epsilon}} \delta$$

uniformly in $\pi$, since (108) clearly follows from the above two statements, for a sufficiently large choice of $C_5$.

We now prove (109). Fix $\pi$, and let $F'$ denote the set

$$F' = \left\{x \in F : \angle(\pi(x), \pi) \lessapprox \delta^{C_5\sqrt{\epsilon}}\right\}.$$

Suppose that $x$ was an element of $F'$. From (104) and the definition of $P_4$ we have

$$\widetilde{\forall}_{y,T} T \in \mathbb{T}[T_\rho], y \in B(x, C\rho) : y \in S(x, \pi(x)).$$

By (31) and (27) we therefore have

$$(110) \qquad \int_{S(x,\pi(x))} \sum_{T \in \mathbb{T}[T_\rho]} \chi_T(y)\, dy \approx \delta^{3/2}.$$

Let $I$ denote the one-dimensional interval in $\mathbb{R}^3$ centered at the origin, parallel to $\pi$, orthogonal to the axis of $T_\rho$, and having length $\delta^{-C_5\sqrt{\epsilon}}\delta$, and let $V$ be a maximal $\delta$-separated subset of $I$. If $Q$ is a square in $S(x, \pi(x))$ and $T$ is a tube in $\mathbb{T}[T_\rho]$, then (66), the angular separation between $\pi(x)$ and $\pi$, and some elementary geometry shows that

$$|T \cap Q| \lesssim \delta^{C_5\sqrt{\epsilon}} \sum_{v \in V} |(T+v) \cap CQ|,$$

where $T + v$ is the translate of $T$ by $v$, and $CQ$ is a dilate of $Q$. Therefore

$$\int_{S(x,\pi(x))} \chi_T(y)\, dy \lesssim \delta^{C_5\sqrt{\epsilon}} \sum_{v \in V} \int_{CS(x,\pi(x))} \chi_{T+v}(y)\, dy$$

where $CS(x, \pi(x))$ was defined in Definition 8.3. Combining this with (110) we obtain

$$\int_{CS(x,\pi(x))} \sum_{v \in V} \sum_{T \in \mathbb{T}[T_\rho]} \chi_{T+v}(y) \, dy \gtrsim \delta^{-C_5\sqrt{\epsilon}} \delta^{3/2}.$$

From (82), (9) and some algebra we therefore obtain

$$\int_{B(x,C\rho)} \left( \sum_{v \in V} \sum_{T \in \mathbb{T}[T_\rho]} \chi_{T+v}(y) \right)^{d'} dy \gtrsim \delta^{C\sqrt{\epsilon}} \delta^{-d'C_5\sqrt{\epsilon}} \delta^{3-d'}.$$

Integrating this over $F'$ we obtain

$$\int |F' \cap B(y, C\rho)| \left( \sum_{v \in V} \sum_{T \in \mathbb{T}[T_\rho]} \chi_{T+v}(y) \right)^{d'} dy \gtrsim \delta^{C\sqrt{\epsilon}} \delta^{-d'C_5\sqrt{\epsilon}} \delta^{3-d'} |F'|.$$

Since $|F' \cap B(y, C\rho)| \lesssim \delta^{3/2}$, we thus obtain

$$(111) \qquad \left\| \sum_{v \in V} \sum_{T \in \mathbb{T}[T_\rho]} \chi_{T+v}(y) \right\|_{d'}^{d'} \gtrsim \delta^{C\sqrt{\epsilon}} \delta^{-d'C_5\sqrt{\epsilon}} \delta^{\frac{3}{2}-d'} |F'|.$$

On the other hand, the directional multiplicity of the collection $\{T+v : T \in \mathbb{T}, v \in V\}$ is at most $\approx \#V \approx \delta^{-C_5\sqrt{\epsilon}}$. Thus by (15), we have

$$\left\| \sum_{v \in V} \sum_{T \in \mathbb{T}[T_\rho]} \chi_{T+v}(y) \right\|_{d'} \lesssim \delta^{\frac{5}{2d'}-1} (\delta^{-C_5\sqrt{\epsilon}})^{1-\beta}$$

for some $\beta > 0$. Combining the two estimates gives (109) with $c = \beta d'$, as desired. □

The above proposition yielded a single arithmetic progression $x_0^0$, $x_1^0$, $x_2^0$ in a single $\rho$-tube $T_\rho$ of length three which had the property $P_5$. By definition of $P_5$, there are many points in $B(x_i^0, C\rho)$, $i = 0, 1, 2$, with the property $P_4$. We therefore expect a large number (about $\delta^{-1}$) of tubes $T$ in $\mathbb{T}[T_\rho]$ to satisfy property $P_4$ when they pass near $x_i^0$.

The property $P_4$ is of course very complicated. However, we shall need only a few of the components of $P_4$ for the purposes of our argument, namely the graininess property, the planiness property at scale $\delta$, and the consistency condition (95). More precisely, the proposition we will need for later sections is as follows.

PROPOSITION 10.2. *Suppose we have an X-ray estimate at dimension $d$, and there exists a Besicovitch set of Minkowski dimension at most $d + \varepsilon$. Then for any $0 < \delta \ll 1$ we can find a $\rho$-tube $T_\rho$, a direction-separated set $\mathbb{T}[T_\rho]$ of* $\delta$-tubes in $T_\rho$ with cardinality $\gtrapprox \delta^{C\sqrt{\epsilon}}\delta^{-1}$, three points $x_i^0$ in arithmetic progression in $T_\rho$ with spacing $\sim 1$, and three planes $\pi_0^0$, $\pi_1^0$, $\pi_2^0$ each containing the axis of $T_\rho$, such that

$$\angle \pi_i^0, \pi_j^0 \gtrapprox \delta^{C\sqrt{\epsilon}} \tag{112}$$

for all $i \neq j$. Furthermore, for each $i = 0, 1, 2$ we can find sets $S_i$ which are grainy with respect to $(x_i^0, \pi_i^0)$, and to each $x_i \in S_i$ we can associate a plane $\pi_i(x_i)$ containing $x_i$, such that

$$\angle \pi_i(x_i), \pi_i^0 \lessapprox \delta^{-C\sqrt{\epsilon}} \rho \tag{113}$$

holds for all $i$ and $x_i \in S_i$, and

$$\int_X \#\{T \in \mathbb{T}[T_\rho] \cap \mathbb{T}(x_0) \cap \mathbb{T}(x_2) : T \in N(\pi_0(x_0)) \cap N(\pi_2(x_2))\} \, dx_0 dx_2 \tag{114}$$
$$\gtrapprox \delta^{C\sqrt{\epsilon}} \delta^4,$$

where

$$N(\pi) = N_{\delta^{-C\sqrt{\epsilon}}\delta}(\pi)$$

and

$$X = \left\{(x_0, x_2) \in S_0 \times S_2 : \frac{1}{2}(x_0 + x_2) \in S_1\right\}.$$

Note that this proposition differs from the previous ones, in that the $\tilde{\forall}$ symbol is completely absent. That stringent notation was necessary only in order to obtain the arithmetic progression in Proposition 10.1. Henceforth we shall be content with properties which are satisfied[1] only about $\delta^{C\sqrt{\epsilon}}$ of the time.

This proposition is not so much a statement on the collection $\mathbb{T}$, but rather on the collection $\mathbb{T}[T_\rho]$, which can be regarded as a rescaled Besicovitch set. We are restricting ourselves to this portion of $\mathbb{T}$ in order to take full advantage of graininess, and will no longer concern ourselves with the remainder of $\mathbb{T}$.

*Proof.* Fix $T_\rho$, $x_0$, $s$ so that the conclusions of Proposition 10.1 hold. Set $x_i^0$ equal to $x_0 + si$ and $S_i$ equal to the set $CS(x_i^0, \pi(x_i^0))$. Set $\pi_i^0$ to be any plane within an angle of $\delta^{-C\sqrt{\epsilon}}\rho$ to $\pi(x_i^0)$ which contains the axis of $T_\rho$; such a plane exists from (66). Note that since $S_i$ is grainy with respect to $(x_i^0, \pi(x_i^0))$, it is also grainy with respect to $(x_i^0, \pi_i^0)$ with slightly worse constants.

Let $i = 0, 1, 2$. From the definition of $P_5(x_i^0, T_\rho)$ we have

$$\tilde{\forall}_{T, x_i} T \in \mathbb{T}[T_\rho], x_i \in T \cap B(x_i^0, C\rho) : P_4(x_i^0, x_i, T, T_\rho).$$

---
[1] See however Proposition 14.5 in the Appendix.

By Lemma 5.6 we therefore have

$$\tilde{\forall}_{T, x_0, x_1, x_2} T \in \mathbb{T}[T_\rho], x_j \in T \cap B(x_j^0, C\rho) \text{ for } j = 0, 1, 2 : P_4(x_i^0, x_i, T, T_\rho).$$

From Lemma 5.2 we thus have

$$\tilde{\forall}_{T, x_0, x_1, x_2} T \in \mathbb{T}[T_\rho], x_j \in T \cap B(x_j^0, C\rho) \text{ for } j = 0, 1, 2 : P_4(x_i^0, x_i, T, T_\rho)$$

holds for all $i$. Applying (31), we obtain

$$\sum_{T \in \mathbb{T}[T_\rho]} \prod_{i=1}^{3} |\{x_i \in T \cap B(x_i^0, C\rho) : P_4(x_i^0, x_i, T, T_\rho)\}| \approx \sum_{T \in \mathbb{T}[T_\rho]} \prod_{i=1}^{3} |T \cap B(x_i^0, C\rho)|.$$

Computing the right-hand side using (28) and (19) we thus obtain

$$(115) \qquad \sum_{T \in \mathbb{T}[T_\rho]} \prod_{i=1}^{3} |\{x_i \in T \cap B(x_i^0, C\rho) : P_4(x_i^0, x_i, T, T_\rho)\}| \gtrapprox \delta^{13/2}.$$

For $i = 0, 2$, define the sets

$$F_i = \{x_i \in S_i : P_3(x_i)\}.$$

For any $x_0, x_2, T$, the contribution to the above sum is $O(\delta^{5/2})$. Furthermore, from elementary geometry and the definition of $P_4$ and $S_i$, we see that this contribution vanishes unless $x_0 \in T \cap F_0$, $x_2 \in T \cap F_2$, and $(x_0, x_2) \in X$. Thus we have

$$\sum_{T \in \mathbb{T}[T_\rho]} |\{(x_0, x_2) \in X \cap ((T \cap F_0) \times (T \cap F_2)) : P_4(x_i^0, x_i, T, T_\rho) \text{ holds for } i = 0, 2\}| \gtrapprox \delta^4.$$

We may rewrite the above as
(116)
$$\int_{X \cap (F_0 \times F_2)} \#\{T \in \mathbb{T}[T_\rho] : P_4(x_0^0, x_0, T, T_\rho), P_4(x_2^0, x_2, T, T_\rho)\} \, dx_0 dx_2 \gtrapprox \delta^4.$$

From the definition of $P_4$, we have

$$\#\{T \in \mathbb{T}[T_\rho] : P_4(x_i^0, x_i, T, T_\rho) \text{ holds for } i = 0, 2\}$$
$$\lessapprox \sum_{\pi_0 \in \Pi_\delta^0(x_0)} \sum_{\pi_2 \in \Pi_\delta^2(x_2)} \#\{T \in \mathbb{T}[T_\rho] \cap \mathbb{T}(x_0) \cap \mathbb{T}(x_2) : T \in N(\pi_0) \cap N(\pi_2)\},$$

where

$$\Pi_\delta^i(x_i) = \{\pi \in \Pi_\delta(x_i) : \angle \pi, \pi_i^0 \lessapprox \delta^{-C\sqrt{\epsilon}}\rho\}.$$

It therefore follows that we may find a way of assigning to each $x_i \in F_i$, $i = 0, 2$ a plane $\pi_i(x_i) \in \Pi_\delta^i(x_i)$ such that

$$(117) \qquad \int_{X \cap (F_0 \times F_2)} \#\{T \in \mathbb{T}[T_\rho] \cap \mathbb{T}(x_0) \cap \mathbb{T}(x_2) :$$
$$T \in N(\pi_0(x_0)) \cap N(\pi_2(x_2))\} \, dx_0 dx_2 \gtrapprox \delta^{C\sqrt{\epsilon}} \delta^4$$

for some large $C$. For, if (117) failed for all choices of assignments $\pi_0(x_0)$ and $\pi_2(x_2)$, then we could average over all such assignments using (91) and contradict (116) if $C$ was large enough.

We thus choose planes $\pi_i(x_i) \in \Pi^i_\delta(x_i)$ so that (117) holds; from (95) we have (113) for all $x_i \in F_i$. We then extend $\pi_i()$ to all of $S_i$ so that (113) holds for all $x_i \in S_i$. It is then an easy matter to check all the properties in the proposition. □

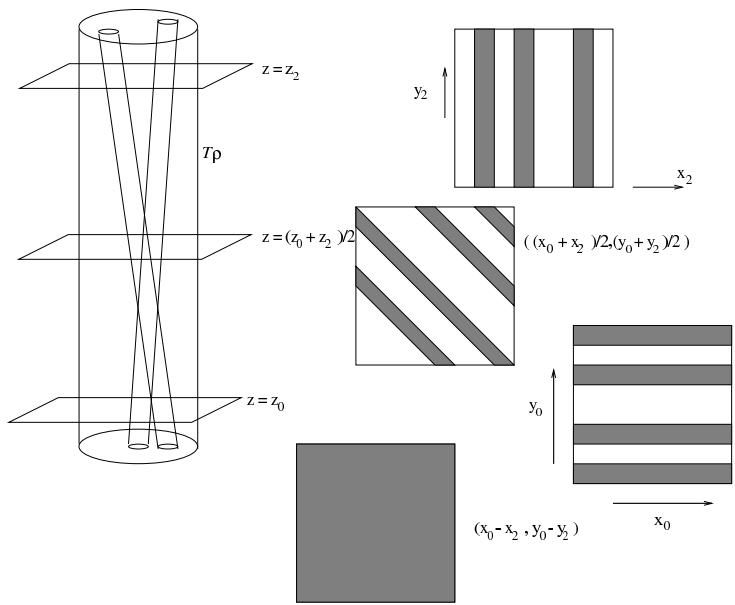

Figure 7. The set $E_\delta[T_\rho]$, and its (grainy) slices at the planes $z = z_0$, $z = z_2$, $z = (z_0 + z_2)/2$. The direction-separated assumption on $\mathbb{T}[T_\rho]$ means that the vectors $(x_0 - x_2, y_0 - y_2)$ must essentially fill out a $\rho \times \rho$ square.

## 11. Reduction to an additive setting

The purpose of this section is to convert the conclusions of Proposition 10.2 into pure additive combinatorics via some affine transformations and a discretization. This follows the philosophy of Bourgain [4], in which the properties of the Besicovitch set were converted into statements about sums and differences of sets; our main innovation is that we shall also encode the planiness and graininess properties into this additive combinatorial setting. After some technical preparations, the desired contradiction will then follow from the lemmas in [4].

For reasons that shall become clear, we shall parametrize the integer lattice $\mathbb{Z}^4$ by $(x_0, y_0, x_2, y_2)$.

PROPOSITION 11.1. *Suppose there exists a Besicovitch set of Minkowski dimension at most $d+\varepsilon$ for some $0 \leq \varepsilon \ll 1$, and let $M \gg 1$ be a large number. Then there exists a subset $\Omega \subset \mathbb{Z}^4$ with cardinality*

$$\#\Omega \gtrsim M^{-C\sqrt{\varepsilon}} M^2 \tag{118}$$

*such that*

- (*Graininess*). *If we define the projections $X_0, Y_0, X_2, Y_2 \in \mathbb{Z}$ of $\Omega$ by*

$$X_0 = \{x_0 : (x_0, y_0, x_2, y_2) \in \Omega \text{ for some } y_0, x_2, y_2\},$$
$$Y_0 = \{y_0 : (x_0, y_0, x_2, y_2) \in \Omega \text{ for some } x_0, x_2, y_2\},$$
$$X_2 = \{x_2 : (x_0, y_0, x_2, y_2) \in \Omega \text{ for some } x_0, y_0, y_2\},$$
$$Y_2 = \{y_2 : (x_0, y_0, x_2, y_2) \in \Omega \text{ for some } x_0, y_0, x_2\},$$

  *then we have*

$$\#X_0 \lesssim M^{C\sqrt{\varepsilon}} M, \tag{119}$$
$$\#Y_0 \lesssim M^{C\sqrt{\varepsilon}} M^{d-2}, \tag{120}$$
$$\#X_2 \lesssim M^{C\sqrt{\varepsilon}} M^{d-2}, \tag{121}$$
$$\#Y_2 \lesssim M^{C\sqrt{\varepsilon}} M. \tag{122}$$

  *Furthermore, we have*

$$\#\{(x_2 + y_0) + (x_0 + y_2) : (x_0, y_0, x_2, y_2) \in \Omega\} \lesssim M^{C\sqrt{\varepsilon}} M^{d-2}. \tag{123}$$

- (*Planiness*). *There exist functions $g : X_0 \times Y_0 \to Y_2$ and $h : X_2 \times Y_2 \to X_0$ such that*

$$y_2 = g(x_0, y_0) \tag{124}$$

  *and*

$$x_0 = h(x_2, y_2) \tag{125}$$

  *for all $(x_0, y_0, x_2, y_2) \in \Omega$.*

- (*Direction-separated tubes*). *We have*

$$\#\{(x_0, y_0, x_2, y_2) \in \Omega : (x_2 + y_0) - (x_0 + y_2) = w\} \lesssim M^{C\sqrt{\varepsilon}} M \tag{126}$$

  *for all $w \in \mathbb{Z}$.*

One should view (126) and (123) as a statement about the difference set and sum set of $X_2 + Y_0$ and $X_0 + Y_2$. The conflict between these two estimates will eventually cause a contradiction if $\varepsilon$ is sufficiently small.

*Proof.* Fix $M$, and define $\delta = M^{-2}$, so that $M = \delta^{-1/2}$. Note that it suffices to prove the above estimates with $\lesssim$ replaced by $\lessapprox$.

For technical reasons we shall prove the proposition with the lattice $\mathbb{Z}$ replaced by $\delta\mathbb{Z}$, and $\mathbb{Z}^4$ by $(\delta\mathbb{Z})^4$. Since the two lattices are isomorphic, this will not make any difference.

By the results of the previous section, we can find $T_\rho$, $\mathbb{T}(T_\rho)$, $x_i^0$, $\pi_i^0$, $S_i$, and $\pi_i()$ which obey the conclusions of Proposition 10.2.

We now apply some affine transformations to place these objects into a standardized form. Firstly we may apply an affine transformation (preserving the z coordinate) so that $T_\rho$ is the vertical tube

$$T_\rho = \{(\mathrm{x},\mathrm{y},\mathrm{z}) : |\mathrm{z}| \lesssim 1, \quad |\mathrm{x}|, |\mathrm{y}| \lesssim \rho\}.$$

By construction and (112), the planes $\pi_i^0$ all contain the z-axis, and make an angle of $\gtrapprox \delta^{C\sqrt{\epsilon}}$ with respect to each other.

By a rotation around the z-axis we may make $\pi_0^0$ the xz-plane:

(127) $$\pi_0^0 = \{(\mathrm{x},\mathrm{y},\mathrm{z}) : \mathrm{y} = 0\}.$$

By (112), we may apply a linear transformation $(\mathrm{x},\mathrm{y},\mathrm{z}) \mapsto (\mathrm{x} + k\mathrm{y}, \mathrm{y}, \mathrm{z})$ for some $|k| \lesssim \delta^{-C\sqrt{\epsilon}}$ and make $\pi_2^0$ the yz-plane without affecting $\pi_0^0$:

(128) $$\pi_2^0 = \{(\mathrm{x},\mathrm{y},\mathrm{z}) : \mathrm{x} = 0\}.$$

This transformation may distort the tube $T_\rho$ and the tubes $T$ in $\mathbb{T}(T_\rho)$ slightly by a factor of $O(\delta^{-C\sqrt{\epsilon}})$, but this will not significantly affect the argument which follows.

Using (112) again, we see that by applying the linear transformation $(\mathrm{x},\mathrm{y},\mathrm{z}) \mapsto (k'\mathrm{x}, k''\mathrm{y}, \mathrm{z})$ for some $\delta^{C\sqrt{\epsilon}} \lesssim |k'|, |k''| \lesssim \delta^{-C\sqrt{\epsilon}}$ we may make $\pi_1^0$ equal to the $\mathrm{x} + \mathrm{y} = 0$ plane without affecting $\pi_0^0$ or $\pi_2^0$:

(129) $$\pi_1^0 = \{(\mathrm{x},\mathrm{y},\mathrm{z}) : \mathrm{x} + \mathrm{y} = 0\}.$$

Again this may distort the tubes $T_\rho$ and $\mathbb{T}(T_\rho)$ slightly, but this will not affect the argument.

For any $x_0, x_2 \in \mathbb{R}^3$ let $\mathbb{T}(x_0, x_2)$ denote the set

$$\mathbb{T}(x_0, x_2) = \{T \in \mathbb{T}[T_\rho] \cap \mathbb{T}(x_0) \cap \mathbb{T}(x_2) : T \subset N(\pi_0(x_0)) \cap N(\pi_2(x_2))\}.$$

From (114) we have

$$\int\int\left[\int_{X(\mathrm{z}_0,\mathrm{z}_2)} \#\mathbb{T}(x_0,x_2) d\mathrm{x}_0 d\mathrm{y}_0 d\mathrm{x}_2 d\mathrm{y}_2\right] d\mathrm{z}_0 d\mathrm{z}_2 \gtrapprox \delta^{C\sqrt{\epsilon}}\delta^4,$$

where $x_i = (\mathrm{x}_i, \mathrm{y}_i, \mathrm{z}_i)$ for $i = 0, 2$ and

$$X(\mathrm{z}_0, \mathrm{z}_2) = \{(\mathrm{x}_0, \mathrm{y}_0, \mathrm{x}_2, \mathrm{y}_2) : ((\mathrm{x}_0, \mathrm{y}_0, \mathrm{z}_0), (\mathrm{x}_2, \mathrm{y}_2, \mathrm{z}_2)) \in X\}.$$

The quantities $z_0$, $z_2$ both range inside intervals of length $O(\rho)$. Thus one can find $z_0, z_2 \in \mathbb{R}$ such that

$$\int_{X(z_0,z_2)} \#\mathbb{T}(x_0, x_2) dx_0 dy_0 dx_2 dy_2 \gtrapprox \delta^{C\sqrt{\epsilon}} \delta^3.$$

Fix $z_0$, $z_2$ so that the above estimate holds. By factoring the $(x_0, y_0)$ and $(x_2, y_2)$ planes into translates of the lattice $(\delta\mathbb{Z})^2$, we see that there must exist some translates $L_0$ and $L_2$ of the lattice $(\delta\mathbb{Z})^2$ such that

$$\sum_{(x_0,y_0,x_2,y_2)\in X(z_0,z_2)\cap(L_0\times L_2)} \#\mathbb{T}(x_0, x_2) \gtrapprox \delta^{-C\sqrt{\epsilon}} \delta^{-1}.$$

By a very mild affine transformation we may make $L_0 = L_2 = (\delta\mathbb{Z})^2$.

Define $\Omega$ to be the set of all quadruplets $(x_0, y_0, x_2, y_2)$ which give a nonzero contribution to the above sum. From the direction-separated nature of $\mathbb{T}[T_\rho]$, each such quadruplet contributes at most $\lessapprox \delta^{-C\sqrt{\epsilon}}$ to the sum, and so we see that $\Omega$ satisfies (118).

Since $(x_0, y_0, z_0)$ lies in $S_0$, and $x_0$, $y_0$ live in $\delta\mathbb{Z}$, the bounds (119) and (120) follow from the graininess of $S_0$ and (127). Similarly the bounds (121) and (122) follow from the graininess of $S_2$ and (128).

Also, since $\left(\frac{x_0+x_2}{2}, \frac{y_0+y_2}{2}, \frac{z_0+z_2}{2}\right)$ lies in $S_0$, the bound (123) follows from (129) and the graininess of $S_1$.

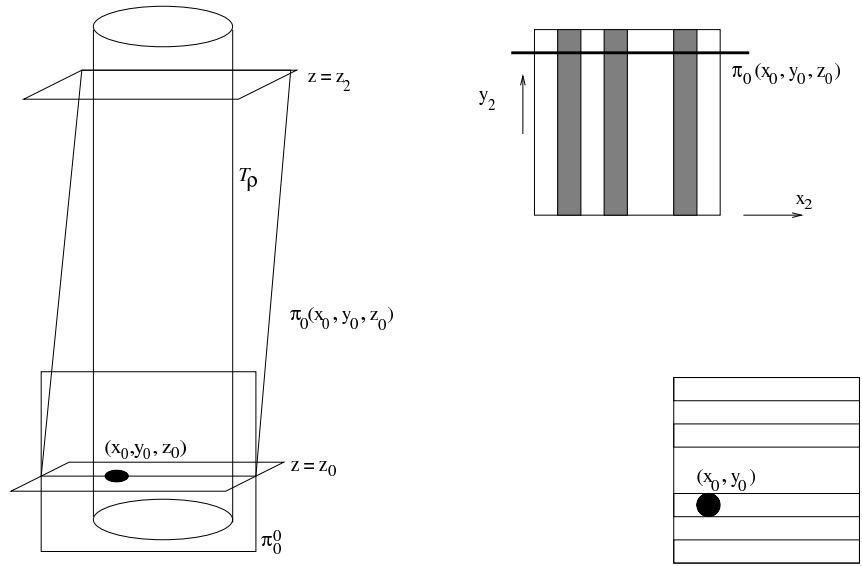

Figure 8. Since the plane $\pi_0(x_0, y_0, z_0)$ must be close to $\pi_0^0$, its intersection with $\{z = z_2\}$ is essentially a horizontal line, which explains (124).

From the direction-separated nature of the tubes in $\mathbb{T}[T_\rho]$ we see that

$$\#\{(x_0, y_0, x_2, y_2) \in \Omega : (x_2 - x_0, y_2 - y_0) = (x, y)\} \lessapprox M^{C\sqrt{\epsilon}}$$

for all $x, y$. The estimate (126) then follows from this and the fact that $x_2 - x_0$ ranges over a $\delta$-separated subset of a $O(\delta^{-C\sqrt{\epsilon}}\delta^{1/2})$-interval, which therefore has cardinality $O(M^{C\sqrt{\epsilon}}M)$.

Finally, we turn to (125) and (124). Fix $(x_0, y_0, x_2, y_2) \in \Omega$, and let $T$ be a tube in $\mathbb{T}(x_0, x_2)$, where $x_i = (x_i, y_i, z_i)$. Then by definition of $\mathbb{T}(x_0, x_2)$ we see that $T$ is contained in a $\delta^{-C\sqrt{\epsilon}}\delta$-neighbourhood of the plane $\pi_0(x_0, y_0, z_0)$. From (113), (127), and elementary geometry we therefore see that $y_2$ is contained in a set of cardinality $O(M^{C\sqrt{\epsilon}})$ which is determined only by $\pi_0(x_0, y_0, z_0)$ (so it is independent of $x_2$ or $y_2$). By the pigeonhole principle we may therefore refine $\Omega$ so that (124) holds for some function $g$, without affecting (118) (other than a worsening of constants). Similarly we can refine $\Omega$ further to obtain (125) for some function $h$. □

## 12. Sums, differences, and functional relationships

In the previous arguments, the parameter $d$ was not fixed to have any particular value. However, the arguments we shall now use work for $d = 5/2$ only, and we shall fix this value accordingly. Of course, this value of $d$ is admissible by Theorem 2.3.

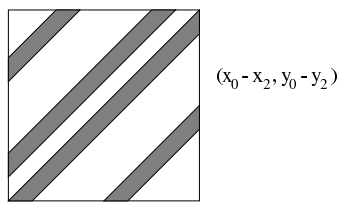

Figure 9. The shape of (130) which is dictated by the all the properties of Proposition 11.1, except for (126), with which it is incompatible.

Having obtained Proposition 11.1, we shall now use pure combinatorics and the lemmas in [4] to show that the conclusions of this proposition lead to a contradiction if $\varepsilon$ is sufficiently small.

The idea is to show that (126) is incompatible with the rest of the properties in Proposition 11.1. In fact, those properties have a strong tendency to make the set

(130) $\quad\{(x_2 - x_0, y_2 - y_0) : (x_0, y_0, x_2, y_2) \in \Omega\}$

resemble the set depicted in Figure 9; this is due to the heuristic implication of (5) from (4), combined with the functional relationships (125) and (124). Of course, such a set is highly incompatible with (126), and so we can attain the desired contradiction.

The property of planiness has evolved into the functional relationships (125) and (124), which state that $x_0$ is independent of $y_0$ and $y_2$ is independent of $x_2$. It turns out that we can strengthen these relationships further in the case $d = 5/2$, and conclude that $y_2$ is also independent of $y_0$, and $x_0$ is independent of $x_2$. In other words, there is a direct functional relationship between $y_2$ and $x_0$.

More precisely:

PROPOSITION 12.1. *Let the notation be as in Proposition 11.1, and suppose that $d = 5/2$. Then there exists a subset $\Omega_1$ of $\Omega$ with cardinality*

$$\#\Omega_1 \gtrsim M^{-C\sqrt{\epsilon}} M^2 \tag{131}$$

*and functions $\tilde{g} : X_0 \to Y_2$ and $\tilde{h} : Y_2 \to X_0$ such that*

$$y_2 = \tilde{g}(x_0) \tag{132}$$

*and*

$$x_0 = \tilde{h}(y_2) \tag{133}$$

*for all $(x_0, y_0, x_2, y_2) \in \Omega_1$. In particular, we have*

$$\#\{(x_2 + y_0) + f(y_2) : (x_0, y_0, x_2, y_2) \in \Omega_1\} \lesssim M^{C\sqrt{\epsilon}} M^{1/2} \tag{134}$$

*and*

$$\#\{(x_0, y_0, x_2, y_2) \in \Omega : (x_2 + y_0) - f(y_2) = w\} \lesssim M^{C\sqrt{\epsilon}} M \tag{135}$$

*for all $w \in \mathbb{Z}$, where $f : Y_2 \to \mathbb{Z}$ is the function $f(y_2) = y_2 + \tilde{h}(y_2)$.*

*Proof.* The argument will be a heavily disguised variant of the "hairbrush" constructions in [17].

For each $(x_0, y_0) \in X_0 \times Y_0$, define the multiplicity

$$\mu_0(x_0, y_0) = \#\{(x_2, y_2) : (x_0, y_0, x_2, y_2) \in \Omega\}.$$

From (121) and (124) we have

$$\|\mu_0\|_\infty \lesssim M^{C\sqrt{\epsilon}} M^{1/2};$$

from (118) we have

$$\|\mu_0\|_1 \gtrsim M^{-C\sqrt{\epsilon}} M^2;$$

and finally, from (119) and (120) we have

$$\#\mathrm{supp}(\mu_0) \lesssim M^{C\sqrt{\epsilon}} M^{3/2}.$$

Combining these facts by Lemma 1.2, we may therefore find a subset $S \subset X_0 \times Y_0$ with cardinality

(136) $$\#S \gtrsim M^{-C\sqrt{\epsilon}}M^{3/2}$$

such that

(137) $$\mu_0(x_0, y_0) \gtrsim M^{-C\sqrt{\epsilon}}M^{1/2} \text{ for all } (x_0, y_0) \in S.$$

In particular, if we define $\Omega' \subset \Omega$ to be the set

$$\{(x_0, y_0, x_2, y_2) \in \Omega : (x_0, y_0) \in S\},$$

then we have

(138) $$\#\Omega' \gtrsim M^{-C\sqrt{\epsilon}}M^2.$$

Now for each $(x_2, y_2) \in X_2 \times Y_2$, define the multiplicity

$$\mu_2(x_2, y_2) = \#\{(x_0, y_0) : (x_0, y_0, x_2, y_2) \in \Omega'\}.$$

From (138) we have

$$\sum_{x_2}\sum_{y_2} \mu_2(x_2, y_2) \gtrsim M^{-C\sqrt{\epsilon}}M^2.$$

From (121), we can therefore find a number $x_2^0 \in X_2$ such that

$$\sum_{y_2} \mu_2(x_2^0, y_2) \gtrsim M^{-C\sqrt{\epsilon}}M^{3/2}.$$

Fix this $x_2^0$. From (120) and (125) we have

$$\sup_{y_2} \mu_2(x_2^0, y_2) \lesssim M^{C\sqrt{\epsilon}}M^{1/2},$$

finally, from (122) we have

$$\#\text{supp}(\mu_2(x_2^0, \cdot)) \lesssim M^{C\sqrt{\epsilon}}M.$$

Combining these facts by Lemma 1.2, we may therefore find a set $Y_2' \subset Y_2$ of cardinality

(139) $$\#Y_2' \gtrsim M^{-C\sqrt{\epsilon}}M$$

such that

(140) $$\mu_2(x_2^0, y_2) \gtrsim M^{-C\sqrt{\epsilon}}M^{1/2} \text{ for all } y_2 \in Y_2'.$$

For each $y_2$ in $Y_2'$, define the set

$$A(y_2) = \{(x_0, y_0) \in Y_0 : (x_0, y_0, x_2^0, y_2) \in \Omega'\}.$$

By (124) these sets are disjoint as $y_2$ varies. By (140) these sets have cardinality

(141) $$\#A(y_2) \gtrsim M^{-C\sqrt{\epsilon}}M^{1/2}.$$

By (125) these sets are contained in the row $\{h(x_2^0, y_2)\} \times Y_0$. Combining these three facts with (120) we thus see that

$$\#\{y_2 \in Y_2' : h(x_2^0, y_2) = x_0\} \lesssim M^{C\sqrt{\epsilon}} \text{ for all } x_0.$$

We may therefore find a subset $Y_2''$ of $Y_2'$ with cardinality

(142) $$\#Y_2'' \gtrsim M^{-C\sqrt{\epsilon}} M$$

such that the map $h(x_2^0, \cdot)$ is injective on $Y_2''$.

We now define the set $\Omega_1$ by

$$\Omega_1 = \left\{ (x_0, y_0, x_2, y_2) : (x_0, y_0) \in \bigcup_{y_2 \in Y_2''} A(y_2) \right\}.$$

By construction we have $A(y_2) \subset S$ for all $y_2 \in Y_2'$. Thus from (137), (141), (142) and the disjointness of the $A(y_2)$ we have

$$\#\Omega_1 \gtrsim M^{-C\sqrt{\epsilon}} M^{1/2} \# \bigcup_{y_2 \in Y_2''} A(y_2) \gtrsim M^{-C\sqrt{\epsilon}} M^{1/2} M^{-C\sqrt{\epsilon}} M M^{-C\sqrt{\epsilon}} M^{1/2}$$

which is (131).

We now show (132). Let $(x_0, y_0, x_2, y_2)$ be an element of $\Omega_1$. Then by construction there exists a $y_2'$ such that $(x_0, y_0, x_2^0, y_2') \in \Omega'$. By (125) we have

$$x_0 = h(x_2, y_2) = h(x_2^0, y_2'),$$

but by (124) we have $y_2 = y_2'$. Thus (133) follows by setting $\tilde{h}(y_2) = h(x_2^0, y_2)$. The claim (132) then follows since $\tilde{h}$ is injective by construction. Finally, the claims (134) and (135) follow from (123) and (126) respectively after applying (133). □

Having strengthened the functional relationship, the next step is to reduce the number of values that $x_2 + y_0$ and $f(y_2)$ can take. More precisely:

PROPOSITION 12.2. *Let the notation be as in the previous propositions. Then one can find a subset $\Omega_2$ of $\Omega_1$ with cardinality*

(143) $$\#\Omega_2 \gtrsim M^{-C\sqrt{\epsilon}} M^2$$

*such that, if $A$ and $B$ denote the sets*

$$A = \{x_2 + y_0 : (x_0, y_0, x_2, y_2) \in \Omega_2\}$$

*and*

$$B = \{f(y_2) : (x_0, y_0, x_2, y_2) \in \Omega_2\},$$

*then*

(144) $$\#A, \#B \lesssim M^{C\sqrt{\epsilon}} M^{1/2}.$$

*Proof.* Define the multiplicity function $\mu_1(a)$ on $\mathbb{Z}$ by
$$\mu_1(a) = \{(x_2, y_0) \in X_2 \times Y_0 : x_2 + y_0 = a\}.$$
We consider the problem of estimating the cardinality of
$$(145) \qquad \left\{(x_0, y_0, x_2, y_2) \in \Omega_1 : \mu_1(x_2 + y_0) \leq M^{-C_7\sqrt{\epsilon}}M^{1/2}\right\}$$
where $C_7$ is a constant to be chosen later.

First consider the contribution of a single $y_2 \in Y_2$ to (145). The quantity $x_0$ is fixed by (133). By (134) there are $\lesssim M^{C\sqrt{\epsilon}}M^{1/2}$ possible values of $x_2 + y_0$ which contribute. Thus the total contribution of $y_2$ to the cardinality of (145) is at most
$$(146) \qquad M^{C\sqrt{\epsilon}}M^{1/2}M^{-C_7\sqrt{\epsilon}}M^{1/2}.$$
By (122) we therefore have
$$\#(145) \lesssim M^{C\sqrt{\epsilon}}M^{-C_7\sqrt{\epsilon}}M^2.$$

Similarly, define the multiplicity function $\mu_2(b)$ on $\mathbb{Z}$ by
$$\mu_2(b) = \{y_2 \in Y_2 : f(y_2) = b\},$$
and consider the cardinality of
$$(147) \qquad \left\{(x_0, y_0, x_2, y_2) \in \Omega_1 : \mu_2(f(y_2)) \leq M^{-C_7\sqrt{\epsilon}}M^{1/2}\right\}.$$
We consider the contribution of a single $(x_2, y_0) \in X_2 \times Y_0$ to (147). By (134) there are $\lesssim M^{C\sqrt{\epsilon}}M^{1/2}$ possible values of $f(y_2)$ which contribute. Since $x_0$ is determined by (133), we see that the contribution of each $(x_2, y_0)$ to the cardinality of (147) is at most (146) as before. From (121) and (120) we thus have
$$\#(147) \lesssim M^{C\sqrt{\epsilon}}M^{-C_7\sqrt{\epsilon}}M^2.$$

Define $\Omega_2$ to be the set $\Omega_1$ with (145) and (147) removed. By (131) we see that (143) holds, if we choose $C_7$ large enough. By construction we have
$$\mu_1(a) \gtrsim M^{-C\sqrt{\epsilon}}M^{1/2} \text{ for all } a \in A$$
and
$$\mu_2(b) \gtrsim M^{-C\sqrt{\epsilon}}M^{1/2} \text{ for all } b \in B.$$
On the other hand, from (121), (120), and (122) we have
$$\|\mu_1\|_1 = \#X_2\#Y_0 \lesssim M^{C\sqrt{\epsilon}}M, \quad \|\mu_2\|_1 = \#Y_2 \lesssim M^{C\sqrt{\epsilon}}M,$$
and (144) follows from Chebyshev's inequality. $\square$

Now we can finally apply Bourgain's lemma. Let $G$ denote the subset of $A \times B$ given by
$$G = \{(x_2 + y_0, f(y_2)) : (x_0, y_0, x_2, y_2) \in \Omega_2\}.$$

From (143) and (135) we have

$$\#\{a - b : (a, b) \in G\} \gtrapprox M^{-C\sqrt{\epsilon}}M.$$

On the other hand, from (134) we have

$$\#\{a + b : (a, b) \in G\} \lessapprox M^{C\sqrt{\epsilon}}M^{1/2}.$$

Finally, the cardinality of $A$ and $B$ is controlled by (144). If we now apply Lemma 2.54 of [4] (or Lemma 14.1 in the Appendix of this paper) with $N = M^{C\sqrt{\epsilon}}M^{1/2}$ we obtain a contradiction, if $\varepsilon$ is sufficiently small. Thus there are no Besicovitch sets of Minkowski dimension $\leq 5/2 + \varepsilon$, if $\varepsilon$ is a sufficiently small absolute constant. This (finally!) completes the proof of Theorem 0.4. □

## 13. Remarks

• A large part of the above argument works for any dimension $d$, providing that one has an X-ray estimate at dimension $d$ to begin with. This raises the possibility that one can obtain a much better result than Theorem 0.4 by a bootstrap argument. At present, however, there are two obstructions to this idea. Firstly, the conclusion of Theorem 0.4 is a Minkowski dimension statement, which is considerably weaker than the X-ray estimate hypothesis. Secondly, the additive combinatorial portion of the argument is only effective at dimension close to $5/2$.

• We now consider the question of what happens when we relax the direction-separated assumption on the collection of tubes $\mathbb{T}$. Our discussion shall be fairly informal.

*Definition* 13.1. A family of $\delta$-tubes $\mathbb{T}$ is said to obey the *Wolff axioms* if:

- For any $\delta \leq \sigma \leq 1$ and any $\sigma$-tube $T_\sigma$, there are at most $\lessapprox \sigma^2 \delta^{-2}$ tubes in $\mathbb{T}$ which are contained in $T_\sigma$.

- For any $\delta \leq \sigma \leq 1$, any $\sigma$-tube $T_\sigma$, and any plane $\pi$ there are at most $\lessapprox \sigma \delta^{-1}$ tubes in $\mathbb{T}$ which are contained in $T_\sigma \cap N_{C\delta}(\pi)$.

Clearly, any direction-separated collection of tubes obeys the Wolff axioms. The Kakeya results in [17] extend easily to these types of collections; for instance, we have

$$\left| \bigcup_{T \in \mathbb{T}} T \right| \gtrapprox \delta^{1/2} \tag{148}$$

whenever $\mathbb{T}$ is a collection of tubes with cardinality $\approx \delta^{-2}$ that satisfies the Wolff axioms. In contrast, the arguments in Bourgain [4] do not say anything concerning this generalized setting in any dimension, because these arguments crucially rely on the direction-separated nature of $\mathbb{T}$. For similar reasons the arguments in our paper also do not adapt to this setting even if one can somehow attain the analogue of stickiness. The estimate (126) in particular is critically dependent on the direction-separation of $\mathbb{T}$.

One may then ask if there is any improvement to (148) if one only assumes that $\mathbb{T}$ satisfies the Wolff axioms, or a similar condition which is weaker than direction-separation. It turns out in fact that there is no such improvement in three dimensions if one replaces the underlying field $\mathbb{R}$ by $\mathbb{C}$, and multiplies all relevant exponents by a factor of 2. Indeed, if we take $E$ to be the unit ball of the Heisenberg group:

$$\left\{(z_1, z_2, z_3) \in \mathbb{C}^3 : \mathrm{Im}(z_3) = \mathrm{Im}(z_1 \overline{z_2}), |z_1|, |z_2|, |z_3| \lesssim 1\right\},$$

then we see that $E$ has Minkowski dimension 5, and for any $a, b \in \mathbb{R}$ and $w \in \mathbb{C}$ with $|a|, |b|, |w| \lesssim 1$ the (complex) line segment

$$\{(z, w + az, z\overline{w} + b) : z \in \mathbb{C}, |z| \lesssim 1\}$$

is contained in $E$. We can then find a collection $\mathbb{T}$ of complex $\delta$-tubes in $N_{C\delta}(E)$ with cardinality $\approx \delta^{-4}$ which satisfy the complex version of the Wolff axioms, and such that

$$\left| \bigcup_{T \in \mathbb{T}} T \right| \sim \delta.$$

Informally, the collection $\mathbb{T}$ is also sticky, plany, and grainy in the spirit of the above argument. By "stickiness" we mean that the tubes $\mathbb{T}$ can be covered by a collection $\mathbb{T}_\rho$ of $\rho$-tubes which are not direction-separated, but satisfy the Wolff axioms and have cardinality $\approx \rho^{-2}$. The planiness and graininess properties are best seen by first verifying them at the origin $(0, 0, 0)$ (in which case the relevant plane is the $z_3 = 0$ plane), and then using the transitive action of the Heisenberg group. One can find a set which obeys all the properties of Proposition 11.1 (with the trivial distinction that the lattice $\mathbb{Z}$ is replaced by the Gaussian integer lattice $\mathbb{Z}[i]$) except for the crucial (126), which is lacking because $\mathbb{T}$ is not direction-separated. Namely, we may take

$$\begin{aligned}
X_0 = Y_2 &= \{a + bi \in \mathbb{Z}[i] : |a + bi| \lesssim M^{1/2}\}, \\
X_2 = Y_0 &= \{a \in \mathbb{Z} : |a| \lesssim M^{1/2}\}, \\
\Omega &= \{(\mathrm{x}_0, \mathrm{y}_0, \mathrm{x}_2, \mathrm{y}_2) \in X_0 \times Y_0 \times X_2 \times Y_2 : \mathrm{x}_0 = \overline{\mathrm{y}_2}\}, \\
g(\mathrm{x}_0, \mathrm{y}_0) &= \overline{\mathrm{x}_0}, \\
h(\mathrm{x}_2, \mathrm{y}_2) &= \overline{\mathrm{y}_2}.
\end{aligned}$$

We remark that the set (130) in this example resembles Figure 9 rather than the set in Figure 7. In particular, (126) fails dramatically.

The counterexamples over $\mathbb{C}$ are based on the fact that $\mathbb{C}$ contains a subfield of half the dimension. It seems natural to ask whether similar objects exist in $\mathbb{R}$; for instance, Erdös [5] has posed the question of whether a Borel subring of $\mathbb{R}$ exists with Hausdorff dimension $1/2$. These types of questions are related to the distance set conjecture of Falconer [6] and are also connected to the Furstenburg problem mentioned in [19]. Two of us will pursue this connection in further detail in [9].

- Our main argument relied on only three slices of the Besicovitch set, and discovered that sets such as $X_2$ and $Y_0$ had very good additive and subtractive properties. One may naturally ask whether one could gain more information by looking at more slices. For instance, one can try to utilize planiness at the slice $z = (z_0 + z_2)/2$; this property was not exploited in the arguments above. This property in fact leads to an interesting observation, that a large portion of the set $E[T_\rho]$ is approximately invariant under a $\delta$-separated set of translations with cardinality $\rho^{-1/2}$ in the $(1, 1, 0)$ direction, and can be used to provide an alternate proof of the contradiction obtained in the arguments above.

If one exploits planiness and graininess on other horizontal slices, one can also discover *multiplicative* properties on $X_2$ and $Y_0$. Very roughly speaking, the sets $X_2$ and $Y_0$ are affinely equivalent to a $\rho$-separated set of diameter $O(1)$ which is an approximate ring, in the sense that it is often closed under addition and multiplication (up to an uncertainty of $\rho$). Thus there is (a very weak and heuristic) converse to the connections between $1/2$-dimensional fields and $5/2$-dimensional Kakeya sets mentioned earlier.

We hope to detail these results in a future paper.

- Define a *Nikodym set* is a subset $E$ of $\mathbb{R}^3$ such that for every $x \in \mathbb{R}^3$ there exists a line $l$ containing $x$ such that $l \cap E$ contains a unit line segment. In analogy with the Kakeya conjecture, it is conjectured that all Nikodym sets have full dimension. At first glance, the above results do not appear to have any bearing on Nikodym sets as the associated collection of $\delta$-tubes is not direction-separated. However this is easily resolved, as one can apply a projective transformation to turn a Nikodym set into a Kakeya set of the same dimension, and so Theorem 0.4 is also valid for Nikodym sets; see [13].

- It seems likely that these techniques could be extended to higher dimensions and give an improvement on (1) for all $n \geq 3$. Indeed, the fact that the codimension of a set with dimension $(n+1)/2$ increases rapidly with $n$ is a favourable sign, especially when compared with the analogue of the graininess property in higher dimensions. Of course, for $n$ sufficiently large ($n > 8$) the estimate (2) is superior to (1) anyway. Two of us have completed the first step of this program in [10], in which the results in [18] are generalized to higher dimensions.

## 14. Appendix: A lemma of Bourgain

In this section we prove a weakened version of Lemma 2.54 of [4], namely:

LEMMA 14.1 ([4]). *Let $N \gg 1$, $0 < \varepsilon \ll 1$ and let $A, B$ be subsets of $\mathbb{Z}$ such that*

$$\#A, \#B \leq N. \tag{149}$$

*Then, if $G$ is a subset of $A \times B$ such that*

$$\#\{a+b : (a,b) \in G\} \leq N, \tag{150}$$

*then we have*

$$\#\{a-b : (a,b) \in G\} \lesssim N^{2-\varepsilon} \tag{151}$$

*for some sufficiently small absolute constant $\varepsilon > 0$.*

In [4] this lemma was proven with the explicit constant $\varepsilon = 1/13$. Our argument is copied from [4] but is slightly simplified as we are not trying to optimize the constants. The argument in [4] is based in turn from the work of Gowers [7]. Recently, two of the authors have managed to improve this to $\varepsilon = 1/6$, see [8]. It is not known what the best value of $\varepsilon$ is, but one must have

$$\varepsilon \leq 2 - \frac{\log(6)}{\log(3)} = 0.36907\ldots;$$

see [8], [12].

*Proof.* Let $\varepsilon$ be a small constant to be chosen later. Suppose for contradiction that we had sets $A$, $B$, $G$ for which (151) failed. Then we clearly have

$$\#G \approx N^2; \tag{152}$$

in this section the $\lesssim$ notation will be used with $\delta$ replaced by $N^{-1}$.

PROPOSITION 14.2. *Let $A$, $B$, $G$ satisfy (152) and the hypotheses of Lemma 14.1. Then there exist $\approx N$ numbers $d$ such that*

$$\#\{(a,b) \in A \times B : a - b = d\} \approx N. \tag{153}$$

A number $d$ which satisfies (153) shall be called an $(A, B)$-*popular difference*.

*Proof.* The upper bound in (153) is clear from (149), as is the upper bound on the number of popular differences. Thus it suffices to show the lower bound in both cases.

Let $X$ denote the set in (150). From (152) we have
$$\|\chi_A * \chi_B\|_{L^1(X)} \geq \#G \gtrapprox N^2.$$
From (150) and (9) with $d=2$ we thus have
$$\|\chi_A * \chi_B\|_2 \gtrapprox N^{3/2}.$$
From the identity
$$\|\chi_A * \chi_B\|_2 = \|\chi_A * \chi_{-B}\|_2;$$
we thus obtain
$$\int (\chi_A * \chi_{-B})^2 \gtrapprox N^3.$$
On the other hand, from (149) we have
$$\int \chi_A * \chi_{-B} \lessapprox N^2;$$
thus
$$\int_{\chi_A * \chi_{-B} \gtrapprox N} (\chi_A * \chi_{-B})^2 \gtrapprox N^3$$
for appropriate choices of constants in the $\gtrapprox$ symbols. On the other hand, from (149) again we have
$$\|\chi_A * \chi_{-B}\|_\infty \lessapprox N.$$
Combining the two estimates we obtain
$$\#\{\chi_A * \chi_{-B} \gtrapprox N\} N^2 \gtrapprox N^3,$$
and the claim follows. $\square$

*Definition* 14.3. If $a, a' \in A$, we say that $a$ and $a'$ $(A,B)$-*communicate* if
$$\#\{b \in B : a-b, a'-b \text{ are both } (A,B)\text{-popular differences}\} \gtrapprox N^{-\sqrt{\epsilon}} N.$$

It is quite plausible that under the conclusions of Proposition 14.2, that many pairs in $A$ communicate. We shall go even further, and show that we can find a refinement of $A$ for which *almost all* pairs communicate.

*Definition* 14.4.  A *refinement* of $A$ is any subset $A'$ of $A$ such that $\#A' \approx \#A$.

PROPOSITION 14.5.  *Let $A$, $B$ be sets satisfying* (149) *and the conclusion of Proposition* 14.2. *Then there exists a refinement $A'$ of $A$ such that*
$$\tilde{\forall}_{a'} a' \in A' : \left[ \tilde{\forall}_a a \in A' : a \text{ and } a' \ (A,B)\text{-communicate} \right].$$

Here the notation $\tilde{\forall}$ is as before, but with $\delta$ replaced by $N^{-1}$.

*Proof.* Let $I \subset A \times B$ be the incidence matrix

$$I = \{(a,b) \in A \times B : a - b \text{ is a } (A,B)\text{-popular difference}\};$$

by hypothesis we have

(154) $$\#I \gtrsim N^2.$$

We define the rows $I_a$ and columns $I^b$ of $I$ by

$$I_a = \{b \in B : (a,b) \in I\}$$

and

$$I^b = \{a \in A : (a,b) \in I\}.$$

From (154) we have

$$\sum_{b \in B} \#I^b \gtrsim N^2.$$

By (149) and Lemma 1.2 one can then find a refinement $B'$ of $B$ such that

(155) $$\#I^b \approx N \text{ for all } b \in B'.$$

Fix $B'$, and let $Y \subset A \times A$ be the set of noncommunicative pairs; i.e.

$$Y = \{(a,a') \in A \times A : \#(I_a \cap I_{a'}) \leq N^{-\sqrt{\epsilon}} N\}.$$

Unfortunately, $Y$ may be quite large. However, for many $b \in B'$ the set $(I^b \times I^b) \cap Y$ is quite small. Indeed, we have

$$\sum_{b \in B'} \#((I^b \times I^b) \cap Y) = \sum_{(a,a') \in Y} \#(I_a \cap I_{a'}) \lesssim N^{-\sqrt{\epsilon}} N (\#A)^2 \lesssim N^{-\sqrt{\epsilon}} N^3$$

by (149). On the other hand, from (155) we see that

(156) $$\sum_{b \in B'} \#(I^b \times I^b) \approx N^2 \#B' \approx N^3.$$

Thus we have

$$\tilde{\forall}_{b,a,a'} b \in B', a \in I^b, a' \in I^b : (a,a') \notin Y.$$

By Corollary 5.9 it follows that

$$\tilde{\forall}_b b \in B' : \left[ \tilde{\forall}_{a,a'} a \in I^b, a' \in I^b : (a,a') \notin Y \right],$$

so that we may find a $b \in B'$ such that

$$\tilde{\forall}_{a,a'} a \in I^b, a' \in I^b : (a,a') \notin Y.$$

Fix this $b$. By Corollary 5.9 again we have

$$\tilde{\forall}_{a'} a' \in I^b : \left[ \tilde{\forall}_a a \in I^b : (a,a') \notin Y \right].$$

The claim then follows by setting $A'$ to be those elements of $I^b$ for which

$$\tilde{\forall}_a a \in I^b : (a, a') \notin Y$$

holds. $\square$

A variant of the above proposition is

PROPOSITION 14.6.  *Let $A$, $B$ be sets satisfying* (149) *and the conclusion of Proposition* 14.2. *Then there exists a refinement $B'$ of $B$ such that for all $b' \in B'$ there are $\approx N$ numbers $a \in A$ for which $a - b'$ is a $(A, B)$-popular difference.*

*Proof.* As in Proposition 14.2, the upper bound on the number of $a$ is trivial from (149).

For each $b' \in B$ let $\mu(b')$ denote the quantity

$$\mu(b') = \#\{a \in A : a - b' \text{ is a } (A, B)\text{-popular difference}\}.$$

Then by hypothesis

$$\|\mu\|_1 \gtrapprox N^2.$$

On the other hand, from (149) we clearly have

$$\mathrm{supp}(\mu), \|\mu\|_\infty \lessapprox N.$$

The claim then follows from Lemma 1.2. $\square$

We now are ready to give the key proposition.

PROPOSITION 14.7.  *Let $A$, $B$, $G$ satisfy* (152) *and the hypotheses of Lemma* 14.1. *Then if $\varepsilon$ is sufficiently small, there exists refinements $A'$ and $B'$ of $A$ and $B$ respectively such that*

(157) $$|G \cap (A' \times B')| \gtrapprox N^2$$

*and*

(158) $$|\{a' - b' : (a', b') \in A' \times B'\}| \lessapprox N^{C\sqrt{\epsilon}} N.$$

*Proof.* Define the multiplicity $\mu_G(a)$ on $A$ by

$$\mu_G(a) = \#\{b \in B : (a, b) \in G\}.$$

Then from (152) we have

$$\|\mu_G\|_1 \gtrapprox N^2$$

Thus by (149) and Lemma 1.2 we can find a refinement $A_1$ of $A$ such that

(159) $$\mu_G(a) \approx N \text{ for all } a \in A_1.$$

Fix this $A_1$. Note that $\#(G \cap (A_1 \times B)) \gtrsim N^2$. By Proposition 14.2 there are $\gtrsim N$ $(A_1, B)$-popular differences, so that by Proposition 14.5 there exists a refinement $A'$ of $A_1$ such that

(160) $$\tilde{\forall}_a a \in A' : a', a \ (A_1, B)\text{-communicate}$$

for all $a' \in A'$.

Fix $A'$. From (159) we have

$$\#(G \cap (A' \times B)) \approx N \# A' \approx N^2.$$

Thus if we define the multiplicity $\mu^G(b)$ on $B$ by

$$\mu^G(b) = \#\{a \in A' : (a, b) \in G\},$$

then

$$\|\mu^G\|_1 \gtrsim N^2.$$

By another application of (149) and Lemma 1.2 we can find a refinement $B_1$ of $B$ such that

(161) $$\mu^G(b) \approx N \text{ for all } b \in B_1.$$

Fix this $B_1$. By Proposition 14.6 there exists a refinement $B'$ of $B_1$ such that for every $b' \in B'$ there are $\approx N$ $a \in A'$ such that $a - b'$ is a $(A', B_1)$-popular difference.

Fix $B'$. Combining this with (160) we conclude (if $\varepsilon$ is sufficiently small)

*Fact* 14.8. For every $a' \in A'$ and $b' \in B'$, there are $\approx N$ elements $a \in A'$ such that $a', a$ $(A, B)$-communicate, and $a - b'$ is a $(A, B)$-popular difference.

We can now finish the proof of Proposition 14.7. Since (157) follows from (161), it suffices to show (158).

Fix $a', b'$. From the above fact and the definition of communicativity, one can find $\gtrsim N^{-C\sqrt{\varepsilon}} N^2$ pairs $(a, b) \in A \times B$ such that $a' - b$, $a - b$, and $a - b'$ are $(A, B)$-popular differences. From the identity

$$a' - b' = (a' - b) - (a - b) + (a - b')$$

and the definition of a popular difference, we thus obtain

$$\#\{(a_1, a_2, a_3, b_1, b_2, b_3) \in A^3 \times B^3 : a' - b' = (a_1 - b_1) - (a_2 - b_2) + (a_3 - b_3)\}$$
$$\gtrsim N^{-C\sqrt{\varepsilon}} N^5.$$

Since there are only $\approx N^6$ sextuplets $(a_1, a_2, a_3, b_1, b_2, b_3)$, (158) follows. □

We can now prove Lemma 14.1. We may refine $G$ so that the map $(a, b) \mapsto a - b$ is injective on $G$ without affecting the hypotheses of the lemma. But this injectivity is clearly incompatible with the conclusions of Proposition 14.7 if $\varepsilon$ is sufficiently small. This is our desired contradiction. □

University of Illinois at Chicago, Chicago, IL
*Current address*: Washington University, St. Louis, MO
*E-mail address*: nets@math.wustl.edu

Princeton University, Princeton, NJ
*E-mail address*: laba@math.princeton.edu

UCLA, Los Angeles, CA
*E-mail address*: tao@math.ucla.edu